\documentclass[opre,nonblindrev]{informs3}

\OneAndAHalfSpacedXI 

\usepackage{endnotes}
\let\footnote=\endnote
\let\enotesize=\normalsize
\def\notesname{Endnotes}%
\def\makeenmark{$^{\theenmark}$}
\def\enoteformat{\rightskip0pt\leftskip0pt\parindent=1.75em
  \leavevmode\llap{\theenmark.\enskip}}

\usepackage{natbib}
 \bibpunct[, ]{(}{)}{,}{a}{}{,}%
 \def\bibfont{\small}%
 \def\bibfont{\small}

\TheoremsNumberedThrough     
\ECRepeatTheorems

\EquationsNumberedThrough    

\usepackage{amsmath,amssymb} 
\usepackage{multirow}
\usepackage{authblk}
\usepackage{cite}
\usepackage{url}
\usepackage{pgf}
\usepackage{color}
\usepackage{todonotes}
\usepackage{booktabs} 
\usepackage{algorithm}
\usepackage[noend]{algorithmic}
\usepackage{pgfplots}
\usepackage{caption}
\usepackage{soul}
\usepackage{bbm}
\usepackage{bm}
\usepackage{pdflscape}

\usepackage[absolute,overlay,showboxes]{textpos}
\usepackage{epstopdf}
\usepackage{subcaption}
\usepackage{changepage}
\usepackage{hyperref}

\usepackage{float}
\usepackage{tikz,pgfplots}
\usetikzlibrary{matrix}
\usepgfplotslibrary{groupplots}
\pgfplotsset{compat=newest}
\usetikzlibrary{calc,positioning,shapes,fit,arrows,shadows}
\tikzstyle{line} = [ draw, -latex']
\usetikzlibrary{shapes,arrows}
\usepgfplotslibrary{units}
\usepackage{url}
\usepackage{color}

\usetikzlibrary{shapes.geometric}
\usetikzlibrary{topaths,calc}

\newtheorem{prop}{Proposition}
\newtheorem{thm}{Theorem}
\newtheorem{lem}{Lemma}
\newtheorem{cor}{Corollary}

\renewcommand{\theARTICLETOP}{}
\renewcommand{\theLRHSecondLine}{}
\renewcommand{\theRRHSecondLine}{}
\renewcommand{\LRHFirstLine}{}
\renewcommand{\RRHFirstLine}{}
\renewcommand{\ECLRHFirstLine}{}
\renewcommand{\ECRRHFirstLine}{}

\begin{document}


\RUNAUTHOR{} 

\RUNTITLE{}

\TITLE{Adaptive Two-stage Stochastic Programming with an Analysis on Capacity Expansion Planning Problem}

\ARTICLEAUTHORS{%
\AUTHOR{Beste Basciftci}
\AFF{Department of Business Analytics, Tippie College of Business, University of Iowa, Iowa City, IA, 52242 \EMAIL{\href{mailto:beste-basciftci@uiowa.edu}{beste-basciftci@uiowa.edu}}}
\AUTHOR{Shabbir Ahmed, Nagi Gebraeel}
\AFF{H. Milton Stewart School of Industrial and Systems Engineering, Georgia Institute of Technology, Atlanta, Georgia, 30332}
}

\ABSTRACT{%
\textbf{\textit{Problem Definition:}} Multi-stage stochastic programming is a well-established framework for sequential decision making under uncertainty by seeking policies that {\color{black}can be dynamically adjusted as uncertainty is realized}. 
Often, e.g., due to contractual constraints, such flexible policies are not desirable, and the decision maker may need to commit to a set of actions for a certain number of periods. Two-stage stochastic programming might be better suited to such settings, where first-stage decisions do not adapt to the uncertainty realized. In this paper, we propose a novel alternative approach, 
named as adaptive two-stage stochastic programming, where each
component of the decision policy requiring limited flexibility has {\color{black}its own revision point, a period prior to
which the decisions are determined at the beginning of the planning until this revision point, and after which they are revised
for adjusting to the uncertainty realized thus far until the end of the planning.
We then analyze this approach} over the capacity expansion planning problem, that may require limited flexibility over expansion decisions.  
\textbf{\textit{Methodology/results:}}
We provide a generic mixed-integer programming formulation for the adaptive two-stage stochastic programming problem {\color{black}with finite support, in particular for scenario trees, }
and show that 
this problem is NP-hard in general. 
Next, we focus on the capacity expansion planning problem, and derive bounds on the value of adaptive two-stage programming in comparison to the two-stage and multi-stage approaches in terms of revision points. 
We propose several heuristic solution algorithms based on this bound analysis. These algorithms either  provide approximation guarantees or computational advantages in solving the resulting adaptive two-stage stochastic 
problem. 
\textbf{\textit{Managerial implications:}} We provide insights on the choice of the revision times based on our analytical analysis. 
We further present an extensive computational study on a generation capacity expansion planning problem with different generation resources including renewable energy. We demonstrate the value of adopting adaptive two-stage approach against the existing policies under limited flexibility and highlight the efficiency of the proposed heuristics along with practical implications on the studied problem.
}%


\maketitle

%

\vspace{-9mm}
\section{Introduction}

Optimization in sequential decision-making processes under uncertainty is known to be a challenging task {\color{black}for addressing multi-period problem settings}. Two-stage and multi-stage stochastic programming are fundamental techniques for modeling these processes, where stage refers to the decision times in planning. In the two-stage approach, the decisions are partitioned into two predetermined portions of the planning horizon - the first stage and second stage. The decisions in the first stage are determined \textit{here-and-now} before the underlying uncertainty {\color{black}is realized}, while those in the second stage are \textit{wait-and-see} decisions that are adjusted with respect to the uncertainty realized.  
In two-stage programs, the first-stage decisions need {\color{black}to be fixed} 
at the beginning of the planning, resulting in \textit{static} policies for those decisions, while the second-stage decisions are the ones that can be adjusted as the {\color{black}uncertainty } information becomes available. On the other hand, multi-stage programs allow total flexibility by deriving \textit{fully adaptive} policies depending on the observed uncertainty at any stage. Thus, both approaches have their own pros and cons. 
Specifically, not allowing any changes for a certain set of decisions while information is revealed over time, such as the first-stage decision structures in two-stage stochastic programming approaches, may lead to costly and rigid solutions.  
However, adopting fully adaptive approaches for all decisions may not be suitable {\color{black}in some practical applications} due to the commitment or lead time related restrictions that the decision makers need to abide by.
In particular, one may need to have a fixed {\color{black}set of decisions before and after a specified time period due to these practical considerations, whose time should be carefully determined. 
To address these issues for such problems with limited flexibilities,} we propose a \textit{partially adaptive} stochastic programming approach that determines the best time to revise the decisions.

\textcolor{black}{There have been many problems in the operations management literature that require partially adaptive policies for determining the best set of actions over a multi-period planning horizon. Supply contracts between buyers and suppliers constitute an important class of problems to determine purchasing and production related decisions over multiple time periods by necessitating a certain set of commitments \citep{Kouvelis2006}. If a buyer desires to update her purchasing commitments considering the uncertainties in the upcoming periods, the supplier might be able to only allow limited flexibility to prevent frequent disruptions in plans. On the other hand, if a buyer and supplier agree on a static policy, this might not be a cost-effective solution over {\color{black} the} long run \citep{Urban2000}. 
\citet{Bassok1997} present the importance of this {\color{black}limited flexibility} in electronics industry by discussing how flexibility levels settled by the buyers and suppliers become critical in determining purchasing policies. Similarly, \citet{BarnesSchuster2002} motivate the necessity of these limitations by providing examples from toys, apparel, and electronics industries. More specifically, the suppliers may not be able to allow multiple updates in purchasing commitments possibly due to production or procurement related lead times, and additional cost of holding unnecessary levels of inventory. 
In these problem settings, it becomes critical to understand how these contracts can be constructed by taking into account the changes at a cost to the buyers and suppliers under this limited flexibility \citep{BenTal2005, Bicer2016}.} 

\textcolor{black}{Partially flexible settings become further necessary in capacity expansion planning problem. This strategic-level planning problem determines primarily the expansion times and {\color{black} expansion} amounts of different resources along with the secondary operational decisions {\color{black}for various application areas including electricity network expansion and production planning} \citep{Luss1982, VanMieghem2003}. These decisions are affected by demand and investment costs that generally involve uncertainties over a multi-period planning horizon.
Setting the expansion decisions at the beginning of the planning results in restrictive and costly  policies {\color{black}as in two-stage models} that consider expansion amounts as the first-stage decisions even if fully flexible policies are allowed for the operational decisions 
\citep{Riis2004, Huang2009}. Generation expansion planning in energy systems is an example problem setting to this case where having static expansion decisions become {\color{black}inconvenient in practice} \citep{Konstantelos2017}. On the other hand, capacity expansion actions might not be updated in each period as in multi-stage models since the expansion of resources may be restricted by commitments in the contracts, or lead time might be a requisite for establishing the necessary infrastructure for traditional generation resources such as nuclear and coal plants that may require more than five years commitment for building them up \citep{Careri2011, Munoz2014}. Therefore, fully adaptive policies obtained from multi-stage models may not be feasible as the expansion decisions of the generation resources need to be set for a certain amount of time. Consequently, these decisions need to take into account the flexibility level of the underlying expansion planning process along with the source of the uncertainties \citep{Latorre2003}.}

\begin{figure}[h]
    \centering
    \caption{{\color{black}Decision dynamics in different methodologies in terms of the duration of the here-and-now decisions determined at the beginning of the planning.}}
    \includegraphics[scale=0.75]{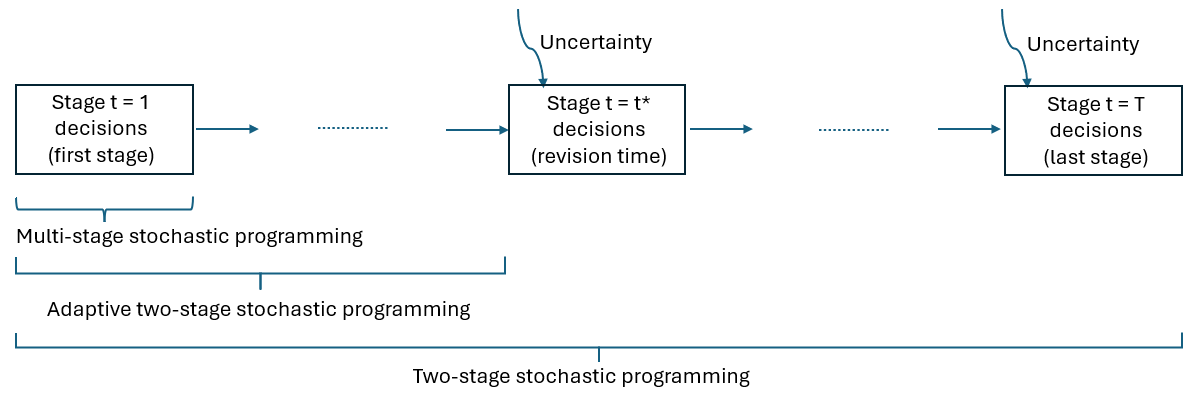}
\label{fig:AdaptiveTSVisualisation}
\end{figure}

These different problem settings motivate the need to incorporate limited  flexibility within the decision-making processes. 
Here, two-stage stochastic programming is a {\it restriction} of a partially adaptive setting by providing a static policy for the first-stage variables, whereas multi-stage stochastic programming is a {\it relaxation} to the desired setting by resulting in a fully flexible solution. 
To address problem settings requiring limited flexibility, we propose a partially adaptive stochastic programming approach, 
in which the stage times are not predetermined but considered as part of the optimization problem. 
{\color{black}Specifically, for each decision with limited adaptability, a {\it revision point} is determined in advance (as a here-and-now decision), where the corresponding decisions are determined at the beginning of the planning until their own revision times, and the remainder of the decisions until the end of the planning horizon are then adjusted at these revision times with respect to the uncertainty realized thus far.}
These revision points are optimally determined for each decision requiring limited flexibility, which is the novelty of our formulation. 
We refer to this approach as {\it adaptive two-stage stochastic programming}. 
{\color{black}Figure \ref{fig:AdaptiveTSVisualisation} provides a conceptual illustration of the proposed framework over decisions with limited flexibility against the existing stochastic programming approaches.}

Since the optimal solution of a multi-stage stochastic program may not provide a feasible {\color{black}implementable} solution for problem settings requiring limited flexibility, this paper does not aim obtaining solutions for the multi-stage model. The paper rather leverages existing approaches, two-stage and multi-stage stochastic programming, to analyze and approximate solutions for the adaptive two-stage stochastic programs. 
{\color{black}To this end, }we make the following contributions:
\begin{enumerate}
\item \textcolor{black}{We introduce {\color{black}the} novel adaptive two-stage stochastic programming approach, in which we optimally determine revision points for each decision under partial flexibility. 
We illustrate this approach using a multi-period newsvendor problem, and provide a policy under the adaptive setting, motivating the importance of the optimization of the revision times. We then develop a mixed-integer linear programming formulation of the proposed approach over generic problems and finite stochastic processes, and prove the NP-hardness of the resulting stochastic program.}
\item \textcolor{black}{We then focus on a specific class of problems, stochastic capacity expansion planning problems, where expansion decisions have limited flexibility over a multi-period planning horizon. 
We provide analytical analyses and insights on the value of the proposed approach compared to two-stage and multi-stage stochastic programming methods with respect to the choice of the revision decisions. 
{\color{black}These} analyses provide upper bounds on the loss against  multi-stage approaches, and lower bounds on the gain against the two-stage approaches depending on the revision times.} 
\item \textcolor{black}{We propose a two-phase procedure for solving the adaptive two-stage programs, where the first phase obtains the revision times and the second phase solves the resulting program given these revision decisions. We present solution algorithms that leverage our bound analyses over the capacity expansion planning problem for obtaining its revision times. We further provide approximation guarantees of these algorithms for solving the expansion planning problem.}
\item \textcolor{black}{We demonstrate the benefits of the adaptive two-stage approach on a generation expansion planning problem. Our computational study illustrates the relative gain of the proposed approach with significant cost reductions, 
compared to two-stage stochastic programming over different scenario tree structures, with solutions that perform similarly well as the policies generated by multi-stage stochastic programming, despite of the limited flexibility. Our results also highlight the significant run time improvements in approximating the desired problem with the help of proposed solution algorithms. We also analyze a sample generation expansion plan with renewable generation resources to further examine the practical implications of optimizing revision decisions.} 
\end{enumerate}

The remainder of the paper is organized as follows: 
In Section \ref{LitReview}, we present the relevant literature. 
{\color{black}
In Section \ref{scProblemFormulation}, we 
formally introduce the adaptive two-stage stochastic programming over a generic setting. 
}
In Section \ref{scAnalyticalAnalysis}, we study the proposed approach on a class of problems encompassing the capacity expansion problem and present analytical results on its performance in comparison to the existing methodologies. In Section \ref{approxAlgSection}, we develop algorithms for solving adaptive two-stage stochastic capacity expansion problems and derive their approximation guarantees by benefiting from our analytical results. In Section \ref{scComputations}, we present an extensive computational study on a generation expansion planning problem in power systems. Section \ref{scConclusion} concludes the paper with final remarks. 
We present all technical proofs in {\color{black}Appendix \ref{adaptiveNewsvendorProof} and Appendix \ref{appendixOmittedProofs}}.

\section{Literature Review}
\label{LitReview}

We first review the relevant studies in the literature that adopt partially adaptive policies in operations management related problems, in particular for inventory and lot sizing problems. 
Partially adaptive policies are motivated in inventory and lot sizing problems for establishing the coordination and synchronization of the supply chain systems over dynamic policies \citep{Silver1988}. \citet{Bookbinder1988} introduce the \textit{static-dynamic} uncertainty strategy for solving a probabilistic lot sizing problem by selecting the replenishment times at the beginning of the planning horizon and determining the corresponding order quantities at these time points. \citet{Tarim2004} extend this concept by developing a mixed-integer linear programming formulation. Variants of this strategy have been studied 
to address the stochastic lot sizing problem under different demand and cost functions, and service level constraints \citep{Zhang2014, Koca2018}, 
depending on the level of flexibility allowed by the suppliers. These approaches mainly consider purchasing and production decisions related to a single product over a restricted number of time periods to quantify the effects of flexible commitments and options by mainly employing dynamic programming techniques or heuristics \citep{Bassok2008, ChenHuang2010}.
However, optimizing the revision times of each decision requiring limited flexibility and the difficulty of handling various constraints have not been addressed in these studies.

In the capacity expansion planning literature, uncertainties over a multi-period planning horizon are addressed with different methods (see \citet{VanMieghem2003} for an extensive survey).
\citet{Sampath1998} study capacity acquisition decisions and their timing to meet customer demand according to the technological breakthroughs by adopting a dynamic programming based solution methodology. 
Similarly, stochastic dynamic programming has been applied to this problem \citep{Wang2017}, despite its disadvantages in incorporating the additional practical constraints \textcolor{black}{along with the curse of dimensionality}. 
Stochastic programming \textcolor{black}{becomes a} fundamental methodology to address these problems by representing the underlying uncertainty through scenarios and \textcolor{black}{considering complicated problem constraints}. 
For instance, \citet{Riis2004} and \citet{Jin2011} 
model these  problems as two-stage stochastic programs by first determining the capacity expansion decisions and then adapting the capacity allocations with respect to the scenarios. On the other hand, \citet{Ahmed2003} and \citet{Singh2009} 
consider these problems as multi-stage stochastic mixed-integer programs and represent uncertainties through scenario trees. Despite this extensive literature, these studies neglect the need for partial flexibility in capacity expansion planning problems.   

In the stochastic programming literature, intermediate approaches between two-stage and multi-stage models have been studied under different problem contexts to mainly address the computational complexity associated with the multi-stage models. As an example, \textit{shrinking-horizon strategy} involves solving two-stage stochastic programs between predetermined time windows. Specifically, \citet{Dempster2000} and \citet{Balasubramanian2014} consider this strategy to obtain an approximation to multi-period planning problems in the oil industry, and multi-product batch plant under demand uncertainty for chemical processes, respectively. 
Shrinking-horizon strategy is also applied to airline revenue management problem in \citet{Chen2010} by proposing heuristics to determine the resolve points under specific assumptions regarding the stochastic process. 
\textcolor{black}{Resolve heuristics are further employed in dynamic programming based stochastic optimization problems such as inventory planning and revenue management \citep{Secomandi2008}.}
Another two-stage approximation to multi-stage models is presented in \citet{Bodur2018} using Linear Decision Rules by limiting the {\color{black}state} decisions to be affine functions of the uncertain parameters{\color{black}, which is then extended by introducing Lagrangian Dual Decision Rules that can handle mixed-integer state variables \citep{Daryalal2024}}.
As an alternative intermediate approach, 
\citet{Zou2018} propose solving a generation capacity expansion planning problem by first considering a multi-stage stochastic program until a predefined stage, and then representing it as a two-stage program. They also develop a \textit{rolling-horizon heuristic} 
to approximate the multi-stage model. 
Another line of research \citep{Dupacova2006} addresses how many stages to have in a multi-stage stochastic program by contamination technique that limits the deviations from the underlying uncertainty distribution, with numerical analyses for the choice of planning horizon and stages for portfolio management \citep{Bertocchi2006}.  
Different to these studies, 
instead of approximately solving multi-stage stochastic programs, we introduce a partially adaptive stochastic programming approach for specifically addressing {\color{black}practical} problem settings under limited flexibility. 

In various problem settings, decisions impact the underlying uncertainty, resulting in complexities for its modelling and effective solution \citep{Jonsbraten1998}. Under this endogenous uncertainty, decisions can change the distribution of the underlying stochasticity  \citep{Hellemo2018, Basciftci2019_Maintenance, Basciftci2019_DDDR} or alternatively they can affect the timing of information acquisition {\color{black}referred to as information discovery} \citep{Goel2006, Vayanos2011, Gupta2014}, which can be related to our paper. 
For capturing the relationship between the underlying stochastic processes and decisions impacting the timing of information discovery, \citet{Goel2006} present a mixed-integer disjunctive programming formulation based on generated scenarios. As solving this problem involves computational challenges, \citet{Vayanos2011} propose a decision rule approximation to ensure a tractable approach, whereas \citet{Gupta2014} discuss alternative decomposition based solution methodologies by relaxing the constraints coupling the scenario sets. \citet{Apap2017} extend these solution methodologies for problems with both endogenous and exogenous uncertainties. 
Stochastic probing \citep{Singla2018, Chugg2019_Probing} can be also considered in this context, which determines the realization of a set of elements over specific problem structures, depending on their probing decisions with a certain cost. 
These studies focus on the information acquisition decisions by changing the information timing reveal for {\color{black}mainly wait-and-see decisions. 
On the other hand, we introduce a generic partially adaptive stochastic programming framework with a novel formulation, where the information discovery is exogenous as the wait-and-see actions are fully flexible, 
whereas decisions that require limited flexibility have revision times, that are optimized at the beginning of the planning, for adjusting these decisions by considering the uncertainty realized thus far.} Additionally, we provide detailed analyses on the choice of revision times on stochastic capacity expansion planning problems and design solution algorithms by leveraging these analyses.

\section{Adaptive Two-stage Stochastic Programming Formulation, Complexity and Value}
\label{scProblemFormulation}

In this section, we 
first propose a generic formulation of the adaptive two-stage approach.
To present the importance of accurately identifying the revision decisions on a motivating setting, we study a multi-period newsvendor problem and give a closed form expression for identifying the {\color{black}ordering} decisions. 
Next, we focus our analyses on finite stochastic processes through utilizing 
scenario trees and provide the resulting generic adaptive two-stage stochastic programming formulation. 
Finally, we show that solving adaptive two-stage stochastic programs is NP-hard.

\subsection{Generic Structure of Adaptive Two-stage Approach}
\label{genericAdaptiveFormSection}

{\color{black}
\subsubsection{Formulation of adaptive two-stage approach}
\label{genericAdaptiveFormSectionFormulation}

To describe a generic formulation for the adaptive two-stage approach, we consider a sequential decision-making problem with $T$ periods, in which we take into account two sets of decisions. The \textit{state} variables $\{x_t (\xi_{[t]})\}_{t=1}^T$ represent the decisions given the data vector $\xi_{[t]}$, namely the data available until period $t$ 
{\color{black}, which } are used for linking decisions of different time periods to each other. The \textit{stage} variables $\{y_t (\xi_{[t]})\}_{t=1}^T$ correspond to the decisions that are local to period $t$. 
{\color{black} We consider the case where the data vector $\xi_{[T]}$ is random and the corresponding stochastic process has a probability distribution $\mathbb{P}$ and support $\Xi$.}
We assume that each {\color{black}$x_t (\xi_{[t]})$ and $y_t(\xi_{[t]})$} 
are vectors of  dimensions $I$ and $J$, respectively. We formalize the adaptive two-stage approach by allowing one revision decision for each state variable throughout the planning horizon,
{\color{black} where the revision decision for the state variable $i$ is defined as $t^*_i$ for $i \in \{1, \cdots, I\}$}. More specifically, the decision maker determines her decisions for the state variable $x_{it} (\xi_{[t]})$ {\color{black} at time $t=1$} until its revision time $t^*_i$ for every $i \in \{1, \cdots, I\}$. Then, she observes the underlying data until that period, namely $\xi_{[t^*_{i}]}$, and adjusts the decisions of the corresponding state variable {\color{black}to be used at that time until the end of the} planning horizon. 
Combining the above, we formulate the adaptive two-stage {\color{black}stochastic} program as follows:
{\color{black}
\begin{subequations} \label{eq:genericAdaptiveFormulation}
\begin{alignat}{2}
\min_{x,\hat{x}, y, t^*} \quad &  \mathbb{E} \left[ f_1 (x_1, y_1)  + \sum_{t = 2}^T f_t (x_t (\xi_{[t]}), y_t (\xi_{[t]}), \xi_{[t]}) \right] \label{eq:genericAdaptiveObj} \\
\text{s.t.} \quad & (x _1, y _1) \in \mathcal{Z}_1, \label{eq:genericAdaptiveRelationship0}\\
& (x_t (\xi_{[t]}), y_t (\xi_{[t]})) \in \mathcal{Z}_t (x _1, \cdots, x_{t-1} (\xi_{[t-1]}),\xi_{[t]}) \quad \mathbb{P}\text{-a.e.} \> \xi_{[T]} \in \Xi, \> t = 2, \cdots, T, \label{eq:genericAdaptiveRelationship} \\
&  x_{it} (\xi_{[t]}) = \hat{x}^{before}_{it} \quad \mathbb{P}\text{-a.e.} \> \xi_{[T]} \in \Xi, \> t = 1, \cdots, t^*_{i} - 1, \> i = 1, \cdots, I, \label{eq:genericAdaptiveConstr1}\\
&  x_{it} (\xi_{[t]}) = \hat{x}^{after}_{it} (\xi_{[t^*_{i}]}) \quad \mathbb{P}\text{-a.e.} \> \xi_{[T]} \in \Xi, \> t = t^*_{i}, \cdots, T, \> i = 1, \cdots, I, \label{eq:genericAdaptiveConstr2} \\
& t_i^* \in \{1, \cdots, T\} \quad \forall i = 1, \cdots, I, \label{eq:genericAdaptiveConstrRevisionNum2}
\end{alignat}
\end{subequations}
}
where 
the function $f_t$ corresponds to the objective function at period $t$.  
{\color{black} Constraint \eqref{eq:genericAdaptiveRelationship0} and the set $\mathcal{Z}_t$ in constraints \eqref{eq:genericAdaptiveRelationship} correspond to the set of constraints at each stage $t =1, \cdots, T$, where $\xi_{[1]}$ is considered as constant as data in the first period is deterministic. 
We introduce the auxiliary decision variables $\hat{x}^{before}$ and $\hat{x}^{after}$ for representing the adaptive relationship depending on the revision time of each state variable $i = 1, \cdots, I$, 
before and after its revision time, as illustrated in constraints \eqref{eq:genericAdaptiveConstr1} and \eqref{eq:genericAdaptiveConstr2}, respectively. Thus, 
underlying data $\xi_{[t^*_i]}$ is observed at stage $t^*_i$, and the corresponding decisions for the remainder of the planning horizon depend on those observations. 
Constraint 
\eqref{eq:genericAdaptiveConstrRevisionNum2} 
represents the revision time of each state variable during the planning horizon.}

{\color{black}
\begin{remark}
We note that although we focus on the case where the state variables have limited flexibility, depending on the studied problem setting, this limitation can be further considered over stage variables as an extension of our formulation. 
Additionally, current formulation considers revision decisions $t^*_i$ for every $i \in \{1, \cdots, I\}$ that are not dependent to the realized uncertainties, 
by considering them as here-and-now decisions. This structure is essential for problem settings, where the decision-makers require the knowledge of the commitment times ahead of the planning. However, if the decision-makers have the flexibility to adjust their revision times depending on the realized uncertainty, then the proposed framework can be extended in a future study by giving further flexibility to the state variables and allowing revision time variables to dynamically change over time while carrying information between different stages. 
\end{remark}
}

{\color{black}
\subsubsection{A motivating example: the newsvendor problem}

\label{scNewsvendor}
To illustrate the importance of the adaptive two-stage approach and optimizing revision times for problems requiring limited flexibility, we evaluate a sample problem with {\color{black}a} single state variable. 
{\color{black} We first present this setting under a given revision time, and then evaluate the resulting policies under different revision timings.}
To this end, we consider a newsvendor problem with $T$ periods, 
where the decision maker determines the order amount in each period $t$, namely $x_t$, while minimizing the total expected cost over the planning horizon. Our goal is to formulate an adaptive two-stage {\color{black}stochastic} program under a given revision time $t^*$, in which we determine the order schedule until 
$t^*$, then observe the underlying uncertainty and determine the {\color{black}ordering decisions for the} remainder of the planning horizon accordingly. 
As discussed earlier, limited flexibility can be required for such a setting, due to contractual restrictions. 
We consider a unit holding cost at the end of each period $t$ as $h_t$, and assume that stockouts are backordered with a cost $b_t$. We incur an ordering cost $c_t$ per unit in each period $t$ and assume that initial inventory at hand is zero. We also assume that the cost structure satisfies the relationship $c_t - b_t \leq c_{t+1} \leq c_t + h_t$ for $t = 1, \cdots, T-1$, because otherwise it might become more profitable to backorder demand or hold inventory. Additionally, we set $b_T \geq c_T$ to avoid backordering at the end of the planning horizon.  
%
Demand in period $t$, denoted by $d_t$ for $t = 1, \cdots, T$, is assumed to be random
{\color{black}, except for $d_1$ which is assumed to be given}.  
By defining the inventory amount at the end of period $t$ as $I_t$ and considering the inventory relationship $I_t = I_{t-1} + x_t - d_t$, we can eliminate the inventory variable using $I_t := \sum_{t' = 1}^t (x_{t'} - d_{t'})$. 
{\color{black}Note that negative inventory corresponds to the amount backordered.} 
Consequently, the dynamic programming formulation of the adaptive {\color{black}two-stage stochastic} newsvendor problem with a {\it given revision point at $t^*$} can be formulated as follows:
\begin{equation} \label{eq:ATS_DPFormulation}
\min_{x_1, \cdots, x_{t^*-1} \in \mathbb{R}_{+}} \left\{\sum_{t = 1}^{t^*-1} \left( c_t x_t + \mathbb{E}[h_t \max\{\sum_{t' = 1}^{t } (x_{t' }- d_{t' }),0\} - b_t \min\{\sum_{t' = 1}^{t} (x_{t'} - d_{t'}),0\}] \right) + \mathbb{E}[Q_{t^*} (\sum_{t = 1}^{t^*-1} (x_t - d_t))] \right\},
\end{equation}
where 
\begin{equation} \label{eq:ATS_DPFormulation2}
Q_{t^*}(s) = \min_{x_{t^*}, \cdots, x_{T} \in \mathbb{R}_{+}} \left\{\sum_{t=t^*}^T \left(c_t x_t + \mathbb{E}[h_t \max\{s + \sum_{t'=t^*}^t (x_{t' } - d_{t' }), 0\} - b_t \min\{s + \sum_{t'=t^*}^t (x_{t' } - d_{t' }), 0\}] \right) \right\}.
\end{equation}

{\color{black}The following result shows that given a revision time, the solution of the newsvendor problem under the adaptive two-stage stochastic programming approach follows an order up to policy.}

\begin{thm} \label{ATSNewsvendorProp} 
Order quantity for the adaptive two-stage problem \eqref{eq:ATS_DPFormulation} can be represented in the following form: 
\begin{align}
\widetilde{F}_{1,t} (X_{1,t}) & = \frac{-c_t + c_{t+1} + b_t}{h_t + b_t}, & \quad t = 1, \cdots, t^* - 1, \\
\widetilde{F}_{t^*,t} (s_{t^*} + X_{t^*,t}) & = \frac{-c_t + c_{t+1} + b_t}{h_t + b_t}, & \quad t = t^*, \cdots, T,
\end{align}
where $X_{i,j}$ is the order up to levels from periods $i$ to $j$, $D_{i,j} = \sum_{t=i}^{j} d_t$, $\widetilde{F}_{i,j}$ is the cumulative distribution function of $D_{i,j}$, $s_{t^*}$ is the inventory level at time $t^{*}$, and $c_{T+1} = 0$. 
\end{thm}

By leveraging the policy in Theorem \ref{ATSNewsvendorProp}, we illustrate the performance of the adaptive two-stage approach under different revision times. 
We consider $T = 5$, and {\color{black}demand in each period $d_t$ to be independently distributed with a truncated normal distribution with mean $\mu_t = 10$, variance $\sigma^2_t = 4$, and a minimum value of 0, }
and cost values are set to $c_t = 5$, $h_t = 2$ for $t = 1, \cdots, 5$  in stationary setting. We let $b_t = h_t$ for the first 4 periods, and set $b_5 = c_5 + 1$ to ensure backordering is costly in the last period. 
We evaluate each revision policy 
{\color{black} under a set of 1000 demand scenarios sampled from this distribution, and compare them with static and dynamic order up to policies under this same set of scenarios. Specifically, the order schedule is determined ahead of the planning in the static setting whereas ordering decisions can be revised in each period by observing the underlying demand in dynamic policies (see \citet{Zipkin2000} for the details of static and fully dynamic policies over the newsvendor problem). In Figure \ref{revisionTimeIllustrationFigNewsvendor}, we present these results by evaluating \eqref{eq:ATS_DPFormulation} under these different policies over this scenario set, where we further illustrate three cases depending on the 
distribution of the demand parameters.}
For the adaptive approach, revising at {\color{black}the} $4^{th}$ period gives the least cost under stationary and increasing demand instances, whereas revising at {\color{black}the} $3^{rd}$ period gives the least cost under decreasing demand case. 
Thus, if revision times are not optimized or ignored for problem settings {\color{black}requiring} limited flexibility, suboptimal solutions can be obtained, which demonstrates the importance of the adaptive two-stage {\color{black}stochastic programming} approach. 
{\color{black}We note that for illustration purposes, the studied problem setting in this section only has a single state variable. However, as multiple decision variables are considered, the problem becomes more challenging and enumeration of potential revision times will no longer be a viable option.} 

\vspace{-3mm}
\begin{figure}[h]
\caption{Objective values under different revision times.}
\label{revisionTimeIllustrationFigNewsvendor}
\centering
\begin{subfigure}{.3\textwidth}
 \centering
\captionsetup{justification=centering}
\begin{tikzpicture}[scale=0.575]
%
%
%
\begin{axis}[
    xlabel={Revision time},
    ylabel={Objective value},
    xmin=1, xmax=5,
    ymin=275, ymax=310,
    xtick={1,2,3,4,5},
    ytick={275,285,295,305},
    legend pos=south east,
    ymajorgrids=true,
    grid style=dashed,
]

\addplot[
    color=black,
    mark=*, dashdotted
    ]
    coordinates {
    (1, 305.12)(2, 305.12)(3, 305.12)(4, 305.12)(5, 305.12)
    }; \label{TS}
\addlegendentry{\small Static}

\addplot[
    color=blue,
    mark=square
    ]
    coordinates {
    (1, 305.12)(2, 295.29)(3, 293.47)(4, 292.81)(5, 297.87)
    }; \label{ATS}
\addlegendentry{\small Partially Adaptive}

\addplot[
    color=red,
    mark=*, dashed
    ]
    coordinates {
    (1, 280.40)(2, 280.40)(3, 280.40)(4, 280.40)(5, 280.40)
    }; \label{MS}
\addlegendentry{\small Fully Dynamic}
\end{axis}
\end{tikzpicture}
\small 
\caption{\footnotesize Stationary demand. \\ \textcolor{white}{empty space}}
\label{EverythingStationary}
\end{subfigure}
\begin{subfigure}{.3\textwidth}
 \centering
\captionsetup{justification=centering}
\begin{tikzpicture}[scale=0.575]
\begin{axis}[
    xlabel={Revision time},
    ylabel={Objective value},
    xmin=1, xmax=5,
    ymin=375, ymax=410,
    xtick={1,2,3,4,5},
    ytick={375, 385, 395, 405},   
    legend pos=south east,
    ymajorgrids=true,
    grid style=dashed,
]

\addplot[
    color=black,
    mark=*, dashdotted
    ]
    coordinates {
    (1, 405.12)(2, 405.12)(3, 405.12)(4, 405.12)(5, 405.12)
    }; \label{TS}
\addlegendentry{\small Static}

\addplot[
    color=blue,
    mark=square
    ]
    coordinates {
    (1, 405.12)(2, 395.29)(3, 393.41)(4, 392.64)(5, 395.43)
    }; \label{ATS}
\addlegendentry{\small Partially Adaptive}

\addplot[
    color=red,
    mark=*, dashed
    ]
    coordinates {
    (1, 380.27)(2, 380.27)(3, 380.27)(4, 380.27)(5, 380.27)
    }; \label{MS}
\addlegendentry{\small Fully Dynamic}
\end{axis}
\end{tikzpicture}
\small
\caption{\footnotesize Increasing demand \\ ($\mu_t = 10 + 2t$).}
\label{CostStationaryDemandIncr}
\end{subfigure}
\begin{subfigure}{.3\textwidth}
 \centering
\captionsetup{justification=centering}
\begin{tikzpicture}[scale=0.575]
\begin{axis}[
    xlabel={Revision time},
    ylabel={Objective value},
    xmin=1, xmax=5,
    ymin=185, ymax=215,
    xtick={1,2,3,4,5},
    ytick={180, 190, 200, 210},   
    legend pos=south east,
    ymajorgrids=true,
    grid style=dashed,
]

\addplot[
    color=black,
    mark=*, dashdotted
    ]
    coordinates {
    (1, 210.96)(2, 210.96)(3, 210.96)(4, 210.96)(5, 210.96)
    }; \label{TS}
\addlegendentry{\small Static}

\addplot[
    color=blue,
    mark=square
    ]
    coordinates {
    (1, 210.96)(2, 199.88)(3, 197.59)(4, 200.76)(5, 210.45)
    }; \label{ATS}
\addlegendentry{\small Partially Adaptive}

\addplot[
    color=red,
    mark=*, dashed
    ]
    coordinates {
    (1, 188.57)(2, 188.57)(3, 188.57)(4, 188.57)(5, 188.57)
    }; \label{MS}
\addlegendentry{\small Fully Dynamic}
\end{axis}
\end{tikzpicture}
\small
\caption{\footnotesize Decreasing demand \\ ($\mu_t = 10 - 2t$).}
\label{CostStationaryDemandDec}
\end{subfigure}
\end{figure}
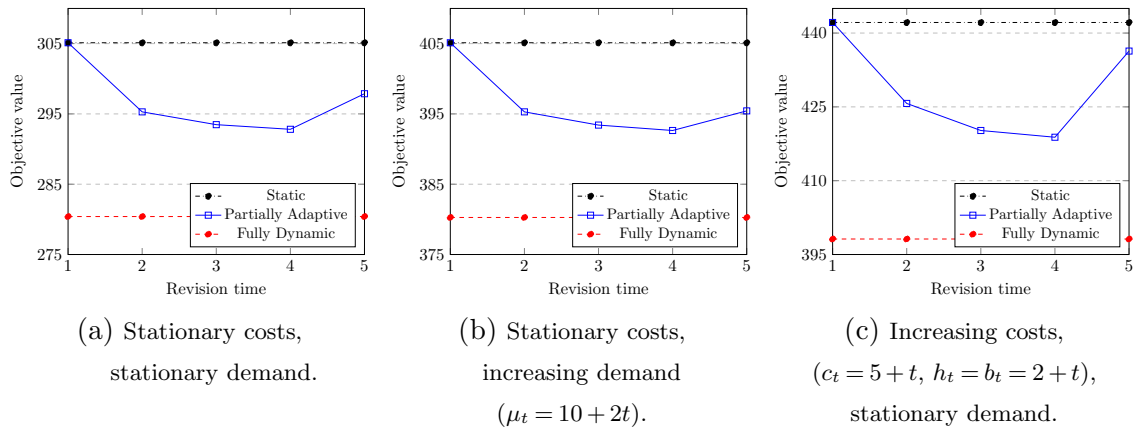 

\subsection{Scenario Tree Formulation} 
In order to represent the adaptive two-stage approach described in Section \ref{genericAdaptiveFormSection} 
{\color{black} over generic problem settings}, we approximate the underlying stochastic process by generating finitely many samples. 
Specifically, a scenario tree is a fundamental method to represent uncertainty in sequential decision-making processes \citep{Ruszczynski2003}, where each node corresponds to a realization of the underlying uncertainty. To construct a scenario tree, sampling approaches similar to the ones proposed for Sample Average Approximation can be adopted, 
where theoretical bounds and confidence intervals regarding the construction of these trees and optimality of the solutions are studied extensively in {\color{black}the} literature (see e.g.,, \citet{Shapiro2003, Kuhn2005}). 

In the remainder of this paper, we consider a scenario tree ${\mathcal{T}}$ with $T$ stages to model the uncertainty structure in a multi-period problem with $T$ periods. 
{\color{black}Note that for notational convenience, the number of stages is considered to be the same as the number of time periods, however these numbers can be different from each other depending on the studied problem setting.}
We illustrate a sample scenario tree in Figure \ref{fig:ScenarioTree}, and represent each node of the tree as $n \in {\mathcal{T}}$. 
We define the set of nodes in each period $1 \leq t \leq T$ as $S_t$, and the period of a node $n$ as $t_n$. Each node $n$, except the root node, has an ancestor, which is denoted as $a(n)$. The unique path from the root node to a specific node $n$ is represented by $P(n)$. We note that each path from the root node to a leaf node corresponds to a scenario, in other words each $P(n)$ gives a scenario when $n \in S_T$. We denote the subtree rooted at node $n$ until period $t$ as ${\mathcal{T}}(n,t)$ for $t_n \leq t \leq T$. To shorten the notation, when the last period of the subtree is $T$, we let ${\mathcal{T}}(n) := {\mathcal{T}}(n,T)$ for all $n \in {\mathcal{T}}$. The probability of each node $n$ is given by $p_n$, where $\sum_{n \in S_t} p_n = 1$ for all $1 \leq t \leq T$. 
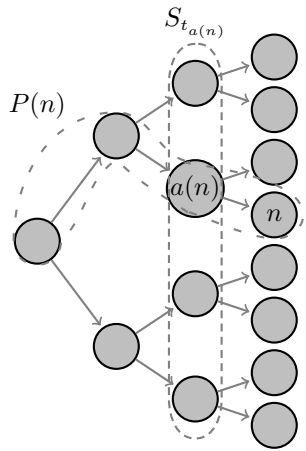
\begin{figure}[H]
\caption{Scenario tree structure.}
\label{fig:ScenarioTree}
\centering
\vspace{-7mm}
\begin{tikzpicture}[shorten >=1pt,->,draw=black!50, node distance=\layersep, thick, scale=0.6]
    \tikzstyle{every pin edge}=[<-,shorten <=1pt]
    \tikzstyle{neuron}=[circle,fill=black!75,minimum size=17pt,draw=black,inner sep=0pt]
    \tikzstyle{input neuron}=[neuron, fill=black!25];
    \tikzstyle{output neuron}=[neuron, fill=red!50];
    \tikzstyle{hidden neuron}=[neuron, fill=blue!50];
    \tikzstyle{annot} = [text width=4em, text centered]

	\node[input neuron] (I0) at (0,0) {};
	\node[input neuron] (I1) at (1.5,2) {};
	\node[input neuron] (I2) at (1.5,-2) {};
	\node[input neuron] (I3) at (3,3) {};
	\node[input neuron] (I4) at (3,1) {\small{$a(n)$}};
	\node[input neuron] (I5) at (3,-1) {};
	\node[input neuron] (I6) at (3,-3) {};
	\node[input neuron] (I7) at (4.5,3.5) {};
	\node[input neuron] (I8) at (4.5,2.5) {};
	\node[input neuron] (I9) at (4.5,1.5) {};
	\node[input neuron] (I10) at (4.5,0.5) {\small{$n$}};
	\node[input neuron] (I11) at (4.5,-0.5) {};
	\node[input neuron] (I12) at (4.5,-1.5) {};
	\node[input neuron] (I13) at (4.5,-2.5) {};
	\node[input neuron] (I14) at (4.5,-3.5) {};
	
            \path (I0) edge (I1);
            \path (I0) edge (I2);
            \path (I1) edge (I3);
            \path (I1) edge (I4);
            \path (I2) edge (I5);
            \path (I2) edge (I6);
            \path (I3) edge (I7);
            \path (I3) edge (I8);
            \path (I4) edge (I9);
            \path (I4) edge (I10);
            \path (I5) edge (I11);
            \path (I5) edge (I12);
            \path (I6) edge (I13);
            \path (I6) edge (I14);

    \begin{scope}[fill opacity=0.0,-, densely dashed]
    \filldraw[fill=yellow!70] ($(I3)+(0,0.75)$) 
        to[out=180,in=90] ($(I3) + (-0.5,0)$) 
        to[out=270,in=90] ($(I6) + (-0.5,0)$)  
        to[out=270,in=180] ($(I6) + (0,-0.75)$)
        to[out=0,in=270] ($(I6) + (0.5,0)$)
        to[out=90,in=270] ($(I3) + (0.5,0)$)  
        to[out=90,in=0] ($(I3) + (0,0.75)$)   ;
        \end{scope}

    \begin{scope}[fill opacity=0.0,-, loosely dashed]
    \filldraw[fill=yellow!70] ($(I0)+(-0.45,0)$) 
    to[out=90,in=180] ($(I1)+(0,0.45)$)
    to[out=0,in=180] ($(I4)+(0,0.5)$) 
    to[out=180,in=90] ($(I10)+(0.55,0)$)
    to[out=270,in=0] ($(I10)+(0,-0.45)$)
    to[out=180,in=145] ($(I4)+(0,-0.55)$)
    to[out=145,in=315] ($(I1)+(0,-0.45)$) 
    to[out=215,in=45] ($(I0)+(0.45,0.0)$) 
    to[out=270,in=0] ($(I0)+(0.0,-0.45)$)
    to[out=180,in=270] ($(I0)+(-0.45,0.0)$) ;
        \end{scope}
        
    \node[annot,above of=I3, node distance=0.8cm] (hl) {\small{$S_{t_{a(n)}}$}};
    \node[annot,above of=I0, node distance=1.8cm] (h) {\small{$P(n)$}};  
\end{tikzpicture}
\end{figure}

We denote the \textit{state variables} as $\{x_n\}_{n \in {\mathcal{T}}}$, and the \textit{stage variables} as $\{y_n\}_{n \in {\mathcal{T}}}$, where \textit{stage variables} $y_n$ are local variables to their associated stage $t_n$. Constraints {\color{black}only} referring to the variables $x_n$ and $y_n$ for each node $n \in {\mathcal{T}}$ are compactly represented by the sets $\mathcal{X}_n$ and $\mathcal{Y}_n${\color{black}, respectively}. We let the dimensions of the variables $\{x_n\}_{n \in {\mathcal{T}}}$ and $\{y_n\}_{n \in {\mathcal{T}}}$ be $I$ and $J$, respectively. We formalize the adaptive two-stage approach in \eqref{eq:adaptiveTwostageFormulationNonlinear} by allowing one revision decision for each state variable throughout the planning horizon. The revision times are denoted by the variable $t^*$. Utilizing the scenario tree structure introduced above, a general form adaptive two-stage stochastic program can be formulated as 
\begin{subequations}  \label{eq:adaptiveTwostageFormulationNonlinear}
\begin{alignat}{1}
\min_{x, y, t^*} \quad &  \sum_{n \in {\mathcal{T}}} p_n (\sum_{i =1}^I a_{in} x_{in} + \sum_{j = 1}^J b_{jn} y_{jn})  \\
\text{s.t.} \quad & \sum_{m \in P(n)} C_{nm} x_m + D_n y_n \geq d_n \quad \forall n \in {\mathcal{T}}, \label{eq:compactConstraint1} \\
& x_n \in \mathcal{X}_n, y_n \in \mathcal{Y}_n, \quad \forall n \in {\mathcal{T}},  \label{eq:compactConstraint2} \\ 
& x_{im} = x_{in} \quad \forall m, n \in S_t, \quad t < t_i^*, i = 1, \cdots, I, \label{eq:constrBeforeRevision1Nonlinear} \\
& x_{im} = x_{in} \quad \forall m, n \in S_t \cap {\mathcal{T}}(j), \quad j \in S_{t_i^*}, \quad t \geq t_i^*, i = 1, \cdots, I,  \label{eq:constrAfterRevision1Nonlinear} \\
& t_i^* \in \{1, \cdots, T\} \quad \forall i = 1, \cdots, I, \label{eq:tStarConstraintNonlinear}
\end{alignat}
\end{subequations}
where the decision variables corresponding to node $n \in {\mathcal{T}}$ are given as $x_n$, $y_n$, and {\color{black}the parameters of node $n\in {\mathcal{T}}$ are} represented as ($a_n$, $b_n$, $C_n$, $D_n$, $d_n$). Constraints \eqref{eq:compactConstraint1} and \eqref{eq:compactConstraint2} describe the constraints over each node $n \in \mathcal{T}$, which are {\color{black}similar to} 
\eqref{eq:genericAdaptiveRelationship0} and \eqref{eq:genericAdaptiveRelationship}. Constraints \eqref{eq:constrBeforeRevision1Nonlinear} and \eqref{eq:constrAfterRevision1Nonlinear} refer to the adaptive two-stage relationship under the revision decisions $t_i^*$, similar to \eqref{eq:genericAdaptiveConstr1} and \eqref{eq:genericAdaptiveConstr2}.  

As we consider problem settings where the limited flexibility is only required for a certain set of decisions, we consider state variables as representatives of these decisions and evaluate different approaches with respect to their adaptability level over these decisions. Two-stage stochastic programming provides a restriction to the adaptive two-stage stochastic programs \eqref{eq:adaptiveTwostageFormulationNonlinear} 
{\color{black}
by adding the constraint \eqref{eq:TwoStageConstr} for each state variable
\begin{equation} \label{eq:TwoStageConstr} 
x_{im} = x_{in} \quad \forall m, n \in S_t, \quad t = 1, \cdots, T,  \quad i = 1, \cdots, I,
\end{equation}
}
which ensures that the state variables $\{x_n\}_{n \in {\mathcal{T}}}$ {\color{black}at a stage} are the same across different scenarios. On the other hand, multi-stage stochastic programming is a relaxation to the adaptive two-stage stochastic programs \eqref{eq:adaptiveTwostageFormulationNonlinear} by {\color{black}removing} the constraints \eqref{eq:constrBeforeRevision1Nonlinear}, \eqref{eq:constrAfterRevision1Nonlinear} and \eqref{eq:tStarConstraintNonlinear} {\color{black}
for all state variables}. 

We illustrate the adaptive two-stage stochastic programming in comparison to existing approaches in Figure \ref{fig:StochProgComparisonFigure} by visualizing the decision structures of each approach over the state variables $\{x_n\}_{n \in {\mathcal{T}}}$, where all approaches assess the same scenario tree for $T$ = 4 periods with two branches from each node (except the leaf nodes) in each stage. To simplify the illustration in terms of the revision point, we consider the state variable $x_n$ as a single variable. 
In the two-stage approach, state variables are set for every time period at the beginning of the planning horizon for each possible realization, as depicted in Figure \ref{fig:TSScnerarioTree}, whereas these variables can be adjusted in multi-stage setting, as shown in Figure \ref{fig:multiStageScnerarioTree}, as uncertainty is revealed over time. 
{\color{black} To illustrate further, the state variables for each node in Figure \ref{fig:multiStageScnerarioTree} can be considered separately as $x_1, x_2, x_3, \cdots, x_{14}, x_{15}$, whereas these variables under the two-stage setting in Figure \ref{fig:TSScnerarioTree} can be considered in the form of $x_1, x_{2-3}, x_{4-7}, x_{8-15}$, as the state variables of the same stage share the same values.} 
In the adaptive two-stage approach, we have a static policy for the state variables before the revision stage. Once we acquire information and adjust our decisions, the remainder trees rooted from each node of the revision stage are compressed to have a single decision for the remaining planning horizon corresponding to that node. This procedure is visualized in Figure \ref{fig:ATSScnerarioTree}, when revision time is selected as the third period. 
{\color{black} Consequently, the corresponding decision variables can be considered as $x_1, x_{2-3}, x_4, x_5, x_6, x_7, x_{8-9}, x_{10-11}, x_{12-13}, x_{14-15}$.}
In particular, the critical point of this approach is optimally determining this revision time with respect to the underlying uncertainty and structure of the studied problem. 
{\color{black} Nevertheless, even under a previously determined revision time, the adaptive two-stage approach provides a more flexible policy than the two-stage approach as the state variables can be adjusted at this time, as opposed to the decision structure in Figure \ref{fig:TSScnerarioTree}.}
We note that we allow $\{y_n\}_{n \in {\mathcal{T}}}$ decisions to have a multi-stage decision structure as in Figure \ref{fig:multiStageScnerarioTree} for all three approaches studied.

\begin{figure}[h]
\caption{Decision structures for $\{x_n\}_{n \in {\mathcal{T}}}$ in different stochastic programming approaches.}
\label{fig:StochProgComparisonFigure}
\centering
\vspace{-3mm}
\begin{subfigure}{4.5cm}
\centering
\begin{tikzpicture}[shorten >=1pt,->,draw=black!50, node distance=\layersep, thick, scale=0.5]
    \tikzstyle{every pin edge}=[<-,shorten <=1pt]
    \tikzstyle{neuron}=[circle,fill=black!75,minimum size=14pt,draw=black,inner sep=0pt]
    \tikzstyle{input neuron}=[neuron, fill=black!25];
    \tikzstyle{output neuron}=[neuron, fill=white!00, draw=white];
    \tikzstyle{hidden neuron}=[neuron, fill=blue!50];
    \tikzstyle{annot} = [text width=4em, text centered]

	\node[output neuron] () at (4.5,5) {};
	\node[output neuron] () at (4.5,-5) {};
	
	\node[input neuron] (I0) at (0,0) {\scriptsize{$x_1$}};
	\node[input neuron] (I1) at (1.5,2) {\scriptsize{$x_2$}};
	\node[input neuron] (I2) at (1.5,-2) {\scriptsize{$x_3$}};
	\node[input neuron] (I3) at (3,3) {\scriptsize{$x_4$}};
	\node[input neuron] (I4) at (3,1) {\scriptsize{$x_5$}};
	\node[input neuron] (I5) at (3,-1) {\scriptsize{$x_6$}};
	\node[input neuron] (I6) at (3,-3) {\scriptsize{$x_7$}};
	\node[input neuron] (I7) at (4.5,3.5) {\scriptsize{$x_8$}};
	\node[input neuron] (I8) at (4.5,2.5) {\scriptsize{$x_9$}};
	\node[input neuron] (I9) at (4.5,1.5) {\scriptsize{$x_{10}$}};
	\node[input neuron] (I10) at (4.5,0.5) {\scriptsize{$x_{11}$}};
	\node[input neuron] (I11) at (4.5,-0.5) {\scriptsize{$x_{12}$}};
	\node[input neuron] (I12) at (4.5,-1.5) {\scriptsize{$x_{13}$}};
	\node[input neuron] (I13) at (4.5,-2.5) {\scriptsize{$x_{14}$}};
	\node[input neuron] (I14) at (4.5,-3.5) {\scriptsize{$x_{15}$}};

            \path (I0) edge (I1);
            \path (I0) edge (I2);
            \path (I1) edge (I3);
            \path (I1) edge (I4);
            \path (I2) edge (I5);
            \path (I2) edge (I6);
            \path (I3) edge (I7);
            \path (I3) edge (I8);
            \path (I4) edge (I9);
            \path (I4) edge (I10);
            \path (I5) edge (I11);
            \path (I5) edge (I12);
            \path (I6) edge (I13);
            \path (I6) edge (I14);
    
\end{tikzpicture}
\caption{Multi-stage.}
\label{fig:multiStageScnerarioTree}
\end{subfigure}
\begin{subfigure}{4.5cm}
\centering
\begin{tikzpicture}[shorten >=1pt,->,draw=black!50, node distance=\layersep, thick, scale=0.5]
    \tikzstyle{every pin edge}=[<-,shorten <=1pt]
    \tikzstyle{neuron}=[circle,fill=black!75,minimum size=10pt,draw=black,inner sep=0pt]
    \tikzstyle{input neuron}=[neuron, fill=black!25];
    \tikzstyle{output neuron}=[neuron, fill=white!00, draw=white];
    \tikzstyle{hidden neuron}=[neuron, fill=blue!50];
    \tikzstyle{annot} = [text width=4em, text centered]

	\node[output neuron] () at (4.5,5) {};
	\node[output neuron] () at (4.5,-5) {};
	
	\node[input neuron] (I0) at (0,0) { };
	\node[input neuron] (I1x2) at (1.5,0) { };
	\node[input neuron] (I3) at (3,3) {};
	\node[input neuron] (I4) at (3,1) {};
	\node[input neuron] (I5) at (3,-1) {};
	\node[input neuron] (I6) at (3,-3) {};
	\node[input neuron] (I7x8) at (4.5,3) {};
	\node[input neuron] (I9x10) at (4.5,1) {};
	\node[input neuron] (I11x12) at (4.5,-1) {};
	\node[input neuron] (I13x14) at (4.5,-3) {};

            \path (I0) edge (I1x2);
            \path (I1x2) edge (I3);
            \path (I1x2) edge (I4);
            \path (I1x2) edge (I5);
            \path (I1x2) edge (I6);
            \path (I3) edge (I7x8);
            \path (I4) edge (I9x10);
            \path (I5) edge (I11x12);
            \path (I6) edge (I13x14);

    \node[annot,below of=I6, node distance=0.7cm] (h) {\footnotesize{revision point}};
            
    \begin{scope}[fill opacity=0.0,-,  dashed]
    \filldraw[fill=yellow!70] ($(I3)+(0,0.0150)$) 
        to[out=90,in=270] ($(I6) + (0,-0.0150)$)   ;
        \end{scope}
    
\end{tikzpicture}
\vspace{0.0cm}
\caption{Adaptive two-stage.}
\label{fig:ATSScnerarioTree}
\end{subfigure}
\begin{subfigure}{4.5cm}
\centering
\begin{tikzpicture}[shorten >=1pt,->,draw=black!50, node distance=\layersep, thick, scale=0.5]
    \tikzstyle{every pin edge}=[<-,shorten <=1pt]
    \tikzstyle{neuron}=[circle,fill=black!75,minimum size=10pt,draw=black,inner sep=0pt]
    \tikzstyle{input neuron}=[neuron, fill=black!25];
    \tikzstyle{output neuron}=[neuron, fill=white!00, draw=white];
    \tikzstyle{hidden neuron}=[neuron, fill=blue!50];
    \tikzstyle{annot} = [text width=4em, text centered]

	\node[output neuron] () at (4.5,5) {};
	\node[output neuron] () at (4.5,-5) {};
	
	\node[input neuron] (I0) at (0,0) { };
	\node[input neuron] (I1x2) at (1.5,0) { };
	\node[input neuron] (I3) at (3,0) {};
	\node[input neuron] (I7x8) at (4.5,0) {};

            \path (I0) edge (I1x2);
            \path (I1x2) edge (I3);
            \path (I3) edge (I7x8);          
\end{tikzpicture}
\caption{Two-stage.}
\label{fig:TSScnerarioTree}
\end{subfigure}
\end{figure}
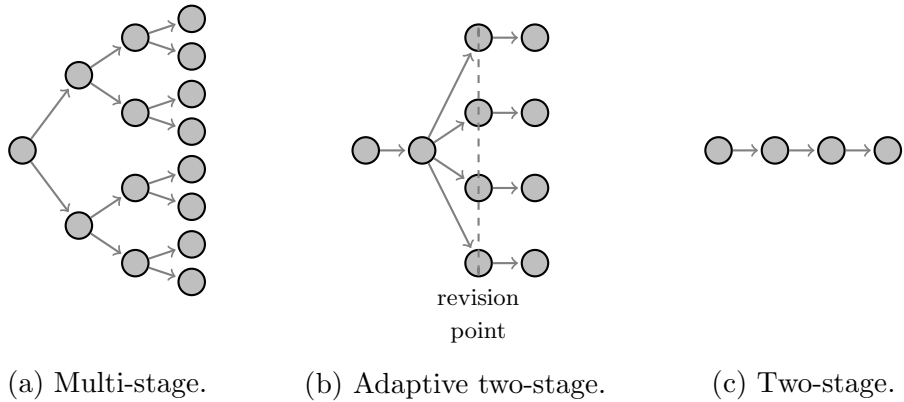


To understand the difficulty of this problem, we identify its computational complexity. 

\begin{thm} \label{thmNPHard}
Solving the adaptive two-stage stochastic programming model in \eqref{eq:adaptiveTwostageFormulationNonlinear} is NP-Hard.
\end{thm}


As constraints \eqref{eq:constrBeforeRevision1Nonlinear} and \eqref{eq:constrAfterRevision1Nonlinear} depend on the decision variable $t_i^*$, we obtain a nonlinear stochastic programming formulation in \eqref{eq:adaptiveTwostageFormulationNonlinear}. To linearize this relationship, we introduce an auxiliary binary variable $r_{it}$ for each $i \in I$, which is 1 if the decisions $\{x_{in}\}$ are revised at nodes $n \in S_t$, and 0 otherwise. Thus, we can reformulate the adaptive two-stage stochastic program as follows:
\begin{subequations}  \label{adaptiveTwostageFormulation}
\begin{alignat}{1}
\min_{x, y, r} \quad &  \sum_{n \in {\mathcal{T}}} p_n (\sum_{i =1}^I a_{in} x_{in} + \sum_{j = 1}^J b_{jn} y_{jn})  \\
\text{s.t.} \quad & \eqref{eq:compactConstraint1}, \eqref{eq:compactConstraint2} \nonumber \\ 
& \sum_{t = 1}^T r_{it} = 1 \quad i = 1, \cdots, I, \label{eq:constrOnceRevision} \\
& x_{im} \geq x_{in} - \bar{x} (1 - \sum_{t' = t+1}^{T} r_{it'}) \quad \forall m, n \in S_t, t = 1, \cdots, T-1, i = 1, \cdots, I, \label{eq:constrBeforeRevision1} \\
& x_{im} \leq x_{in} + \bar{x} (1 - \sum_{t' = t+1}^{T} r_{it'}) \quad \forall m, n \in S_t, t = 1, \cdots, T-1, i = 1, \cdots, I, \label{eq:constrBeforeRevision2} \\
& x_{im} \geq x_{in} - \bar{x} (1 - r_{it}) \quad \forall m, n \in S_{t'} \cap {\mathcal{T}}(l), l \in S_t, t' \geq t, t = 1, \cdots, T, i = 1, \cdots, I, \label{eq:constrAfterRevision1} \\
& x_{im} \leq x_{in} + \bar{x} (1 - r_{it}) \quad \forall m, n \in S_{t'} \cap {\mathcal{T}}(l), l \in S_t, t' \geq t, t = 1, \cdots, T, i = 1, \cdots, I,  \label{eq:constrAfterRevision2} \\
& r_{it} \in \{0,1\} \quad \forall i = 1, \cdots, I, t = 1, \cdots, T, 
\end{alignat}
\end{subequations}

\noindent {\color{black}where $\bar{x}$ represents an upper bound on the values of the state variables. 
Here, constraint \eqref{eq:constrOnceRevision} ensures that each state variable is revised once within the planning horizon. Constraints \eqref{eq:constrBeforeRevision1} and \eqref{eq:constrBeforeRevision2} represent the decision structure of the state variables before their revision times by linearizing the relationship presented in constraint \eqref{eq:constrBeforeRevision1Nonlinear}, whereas constraints \eqref{eq:constrAfterRevision1} and \eqref{eq:constrAfterRevision2} correspond to the decision structures after the revision times by building a similar reformulation for constraint~\eqref{eq:constrAfterRevision1Nonlinear}.}

\begin{remark} \label{remarkATSModelSize}
Although optimization model in (\ref{adaptiveTwostageFormulation}) provides a stochastic mixed-integer linear programming formulation for the adaptive two-stage approach, it requires the addition of exponentially many linear constraints in terms of the number of periods for representing the desired relationship. 
Additionally, constraints \eqref{eq:constrBeforeRevision1}--\eqref{eq:constrAfterRevision2} involve big-$M$ coefficients which may weaken its linear programming relaxation. 
Consequently, it becomes computationally challenging to directly solve the formulation (\ref{adaptiveTwostageFormulation}) when the size of the scenario tree becomes larger. 
\end{remark}


\subsection{Value of Adaptive Two-stage Stochastic Solutions}
\label{sec:ValueofATS}

This section presents metrics to assess the performance of the proposed adaptive two-stage approach 
under a given revision decision $t^* \in \mathbb{Z}_{+}^I$. 
To this end, we have the following relationship for the vector $t^* \in \{1, \cdots, T\}^I$ 
\[
V^{MS} \leq V^{ATS}(t^*) \leq V^{TS},
\]
where $V^{MS}$, $V^{ATS}(t^*)$, and $V^{TS}$ correspond to the objective values of the multi-stage, adaptive two-stage under given revision decision $t^*$, and two-stage stochastic programming models, respectively. This relation holds as the two-stage program provides a feasible solution for the adaptive two-stage program under any $t^*$ vector, and the solution of the adaptive two-stage program under any $t^*$ vector is feasible for the multi-stage program. 

To evaluate the performance of the adaptive two-stage approach, we analyze 
$V^{MS} - V^{ATS}(t^*)$ and $V^{TS} - V^{ATS}(t^*)$ for a given $t^*$ vector. 
We refer $V^{TS} - V^{ATS}(t^*)$ as the \textit{value of adaptive two-stage} (VATS) in the remainder of this paper. 
Although multi-stage approach does not provide a feasible solution for the adaptive two-stage approach, we aim at analyzing $V^{MS} - V^{ATS}(t^*)$ as well, to quantify the impact of the revision decisions in adaptive two-stage setting. 

Due to the complexity of the generic adaptive two-stage stochastic problems as discussed in Theorem \ref{thmNPHard} and  Remark \ref{remarkATSModelSize}, in the remainder of this paper, we focus on analyzing the proposed approach over a specific problem structure that requires limited flexibility on a certain set of decisions, where we can obtain bounds on $V^{MS} - V^{ATS}(t^*)$ and $V^{TS} - V^{ATS}(t^*)$ for a given $t^*$ vector to assess the value of the adaptive two-stage approach based on revision decision $t^*$. 

\section{Analysis of Capacity Expansion Planning Problem}
\label{scAnalyticalAnalysis}

In this section and in the remainder of this paper, we study stochastic capacity expansion planning problem that requires limited flexibility on expansion decisions over a multi-period planning horizon.
Our goal is to  derive analytical bounds to evaluate the value of the adaptive two-stage stochastic programming approach using the metrics defined in Section \ref{sec:ValueofATS}. 
We first present the capacity expansion planning problem as a multi-stage stochastic program. Then, we consider the setting requiring limited flexibility through the adaptive two-stage approach under given revision times, and provide an analytical analysis with insights depending on the choice of these times. 
%

\subsection{Problem Formulation}

The capacity expansion planning problem determines the capacity acquisition decisions of the set of resources~$\mathcal{I}$ by allocating the corresponding capacities to the tasks~$\mathcal{J}$, while satisfying the demand of items~$\mathcal{K}$ and capacity constraints. A multi-stage stochastic capacity expansion planning problem with $T$ periods can be written as follows:
\begin{subequations} \label{eq:CapacityExpansionMultiStage}
\begin{alignat}{1}
\min_{x,y} \quad &  \sum_{n \in {\mathcal{T}}} p_n (a_n^\top x_n + b_n^\top y_n) \label{eq:CapacityExpansionMultiStageObj}\\
\text{s.t.} \quad & A_n y_n \leq \sum_{m \in P(n)} x_{m} \quad \forall n \in {\mathcal{T}}, \label{eq:CapacityExpansionMultiStageConstr1} \\
& B_n y_n \geq d_n \quad \forall n \in {\mathcal{T}}, \label{eq:CapacityExpansionMultiStageConstr2}\\
& x_{n} \in \mathbb{Z}_{+}^{|\mathcal{I}|}, y_n \in \mathbb{R}_{+}^{|\mathcal{J}|} \quad \forall n \in {\mathcal{T}},
\end{alignat}
\end{subequations}
where the decision variables $x_n$, $y_n$, and the parameter $p_n$ represent the capacity expansion/acquisition decisions, capacity allocation decisions, and probability corresponding to the node $n$, respectively. By adopting the scenario tree structure, we consider the problem parameters as $(a_n, b_n, A_n, B_n, d_n)$ corresponding to the node $n$\textcolor{black}{, where $a_n \in \mathbb{R}_{+}^{|\mathcal{I}|}$, $b_n \in \mathbb{R}_{+}^{|\mathcal{J}|}$ represent the cost parameters, and $d_n \in \mathbb{R}_{+}^{|\mathcal{K}|}$ corresponds to the demand for every node $n \in \mathcal{T}$. We let $I = |\mathcal{I}|$}. 
Objective \eqref{eq:CapacityExpansionMultiStageObj} minimizes the total cost by considering capacity acquisition and allocation decisions. Constraint \eqref{eq:CapacityExpansionMultiStageConstr1} guarantees that the assigned capacity in each period is at most the available capacity, and constraint \eqref{eq:CapacityExpansionMultiStageConstr2} ensures that the demand is satisfied.  
We note that the problem \eqref{eq:CapacityExpansionMultiStage} can be reformulated as in \citet{Huang2009} by considering its specific substructure 
as follows:
\begin{subequations} \label{eq:CapacityExpansionMultiStageReformulation1}
\begin{alignat}{1}
\min_{y} \quad &  \sum_{n \in {\mathcal{T}}} p_n  b_n^\top y_n + \sum_{i \in \mathcal{I}} Q_i(y) \\
\text{s.t.} \quad & \eqref{eq:CapacityExpansionMultiStageConstr2} \nonumber \\
& y_n \in \mathbb{R}_{+}^{|\mathcal{J}|} \quad \forall n \in {\mathcal{T}},
\end{alignat}
\end{subequations}
where, for each $i \in \mathcal{I}$, and describing $i^{th}$ component of a vector through the notation $[\cdot]_i$,
\begin{subequations} \label{eq:CapacityExpansionMultiStageReformulation2}
\begin{alignat}{1}
Q_i(y) = \min_{x_i} \quad &  \sum_{n \in {\mathcal{T}}} p_n  a_{in} x_{in},  \\
\text{s.t.} \quad &  \sum_{m \in P(n)} x_{im} \geq [A_n y_n]_i \quad \forall n \in {\mathcal{T}}, \\
& x_{in} \in \mathbb{Z}_{+} \quad \forall n \in {\mathcal{T}}.
\end{alignat}
\end{subequations}

In this reformulation, the problem \eqref{eq:CapacityExpansionMultiStageReformulation1} determines the capacity allocation decisions, whereas the problem \eqref{eq:CapacityExpansionMultiStageReformulation2} considers the capacity acquisition decisions given the allocation decisions, which is further decomposed with respect to each {\color{black}resource} item $i \in \mathcal{I}$. To this end, under a given sequence of allocation decisions $\{y_n\}_{n \in {\mathcal{T}}}$, we can obtain the optimal capacity acquisition decisions corresponding to each item $i$ by solving \eqref{eq:CapacityExpansionMultiStageReformulation2}. This problem can be restated by letting $\delta_{in} = [A_n y_n]_i$ 
as follows:
\begin{subequations} \label{eq:CapacityExpansionMultiStageReformulation3}
\begin{alignat}{1}
\min_{x_i} \quad &  \sum_{n \in {\mathcal{T}}} p_n  a_{in} x_{in}  \\
\text{s.t.} \quad & \sum_{m \in P(n)} x_{im} \geq \delta_{in} \quad \forall n \in {\mathcal{T}} \label{eq:SubstructureConstr1} \\
& x_{in} \in \mathbb{Z}_{+} \quad \forall n \in {\mathcal{T}}. \label{eq:SubstructureConstr2}
\end{alignat}
\end{subequations}

We refer to the problem \eqref{eq:CapacityExpansionMultiStageReformulation3} as the {\it single-resource problem}. 
We note that the feasible region of the formulation \eqref{eq:CapacityExpansionMultiStageReformulation3} gives an integral polyhedron under integer $\{\delta_{in}\}_{n \in {\mathcal{T}}}$ values \citep{Huang2009}. 
For comparison, the two-stage stochastic programming formulation of the single-resource problem \eqref{eq:CapacityExpansionMultiStageReformulation3} can be represented as follows:
{\color{black}
\begin{equation}
\label{eq:CapacityExpansionTwoStageReformulation}
\min_{x_i} \{\sum_{n \in {\mathcal{T}}} p_n  a_{in} x_{in}: \eqref{eq:SubstructureConstr1}, \eqref{eq:SubstructureConstr2}, \> x_{im} = x_{in} \> \forall m, n \in S_t, {\color{black}\> t = 1, \cdots, T}\}.
\end{equation}
}
To study this problem under a partially flexible setting, we introduce {\color{black}the adaptive two-stage stochastic programming formulation of the single-resource problem under a {\it given revision point $t_i^*$},} where the capacity acquisition decision of each item $i \in \mathcal{I}$ is determined at the beginning of planning until stage $t_i^* - 1$. Then, the underlying scenario tree at time $t_i^*$ is observed, and the acquisition decisions until the end of the planning horizon are determined. 
The resulting problem can be formulated as follows under a given $t_i^*$ value:
{\color{black}
\begin{equation}
\label{eq:CapacityExpansionAdaptiveReformulation}
\min_{x_i} \{\sum_{n \in {\mathcal{T}}} p_n  a_{in} x_{in}: \eqref{eq:SubstructureConstr1}, \eqref{eq:SubstructureConstr2}, x_{im} = x_{in} \> \forall m, n \in S_t, \> t < t_i^*, \> x_{im} = x_{in} \> \forall m, n \in S_t \cap {\mathcal{T}}(j), \> j \in S_{t_i^*}, \> t \geq t_i^*\}.
\end{equation}
}
We can equivalently represent \eqref{eq:CapacityExpansionAdaptiveReformulation} by condensing the underlying scenario tree based on the $t_i^*$ value. 
Specifically, let the condensed version of tree ${\mathcal{T}}$ for item $i$ be ${\hat{\mathcal{T}}}_i (t_i^*)$, where we denote the set of nodes that are condensed to node $n \in {\hat{\mathcal{T}}}_i (t_i^*)$ as $\hat C_n \subset {\mathcal{T}}$. The set 
$\hat C_n$ is obtained by combining the nodes with the same decision structure as defined by {\color{black} the constraints of the problem \eqref{eq:CapacityExpansionAdaptiveReformulation}}. 
{\color{black}For example, for the decision structure illustrated in Figure \ref{fig:ATSScnerarioTree}, we can consider the following condensed sets by covering the set of nodes from 1 to 15 from the original tree ${\mathcal{T}}$ through 10 condensed nodes from the set ${\hat{\mathcal{T}}}_i (t_i^*)$: $\hat C_1 = \{1\}$, $\hat C_2 = \{2,3\}$, $\hat C_{3} = \{4\}$, $\hat C_{4} = \{5\}$, $\hat C_{5} = \{6\}$, $\hat C_{6} = \{7\}$, $\hat C_7 = \{8,9\}$, $\hat C_8 = \{10,11\}$, $\hat C_9 = \{12,13\}$, $\hat C_{10} = \{14,15\}$.}
Then, the reformulation of \eqref{eq:CapacityExpansionAdaptiveReformulation} is
{\color{black}
\begin{equation}
\label{eq:CapacityExpansionAdaptiveReformulationCondensed}
\min_{x_i} \{\sum_{n \in {\hat{\mathcal{T}}}_i (t_i^*)} \hat p_n  \hat a_{in} x_{in}: \sum_{m \in \hat{P}(n)} x_{im} \geq \hat \delta_{in} \> \forall n \in {\hat{\mathcal{T}}}_i (t_i^*), \> x_{in} \in \mathbb{Z}_{+} \quad \forall n \in {\hat{\mathcal{T}}}_i (t_i^*)\},
\end{equation}
}
where {\color{black}$\hat p_n = 1$}, 
$\hat a_{in} = \sum_{m \in \hat C_n} p_m a_{im}$, $\hat \delta_{in} = \max_{m \in \hat C_n} \{\delta_{im}\}$, {\color{black}and $\hat{P}(n)$ represents the path of node $n \in {\hat{\mathcal{T}}}_i (t_i^*)$ from the root node under the condensed tree}.

\begin{prop} \label{TUAdaptiveProp}
The coefficient matrix of the formulation \eqref{eq:CapacityExpansionAdaptiveReformulationCondensed} is totally unimodular.
\end{prop}
This result helps in the subsequent section to extend the value of the adaptive two-stage approach from the subproblem \eqref{eq:CapacityExpansionAdaptiveReformulation}, or equivalently \eqref{eq:CapacityExpansionAdaptiveReformulationCondensed}, to the capacity expansion planning problem \eqref{eq:CapacityExpansionMultiStageReformulation1} {\color{black}by considering its linear programming relaxation}.

\subsection{Deriving VATS for Capacity Expansion Planning}

{\color{black} In this section, we first analyze the single-resource problem to derive the VATS over this problem depending on the revision time. We then provide insights by utilizing these analytical results by considering its sensitivity to the different problem parameters. Finally, we extend these results to the capacity expansion planning problem and obtain the value of the proposed approach against two-stage and multi-stage stochastic programming policies.}  


\subsubsection{VATS for the single-resource problem}
\label{VATSSingleItem}

First, we derive the VATS for the single-resource problem \eqref{eq:CapacityExpansionMultiStageReformulation3} using its linear programming relaxations under two-stage, adaptive two-stage and multi-stage models. {\color{black} Under given $\{p_n, a_n, \delta_n\}_{n \in \mathcal{T}}$ parameters, we} let $v^M$, $v^T$, $v^R (t^*)$ be the optimal value of the linear programming relaxations of the formulations \eqref{eq:CapacityExpansionMultiStageReformulation3}, \eqref{eq:CapacityExpansionTwoStageReformulation}, and \eqref{eq:CapacityExpansionAdaptiveReformulation}, respectively. 
{\color{black} We note that since $\delta_n$ values are given, its dependency to the capacity allocation decisions $\{y_n\}_{n \in \tau}$ are omitted. Additionally, }
as we focus on the single-resource problem, we omit the resource index $i$ in this section for brevity. 

To construct the VATS, we examine the problem parameters with respect to the underlying scenario tree. Specifically, we define the minimum and maximum costs over the scenario tree as $a_{*} =  \min_{n \in {\mathcal{T}}} \{ a_{n} \}$ and $a^* =  \max_{n \in {\mathcal{T}}} \{ a_{n} \}$, respectively. 
\textcolor{black}{For the analysis of the single-resource problem, we refer $\delta_n$ as the demand parameter, as under such a perspective this problem is equivalent to the stochastic lot sizing problem without a fixed order cost \citep{Guan2006_Formulations}. Thus, our insights on the single-resource problem can relate to the stochastic lot sizing problem as well.}
We denote the maximum demand over the scenario tree as $\delta^* = \max_{n \in {\mathcal{T}}} \{\delta_{n}\}$, and the expected maximum demand over scenarios as $\bar \delta = \sum_{n \in S_T} p_n \max_{m \in P(n)} \{\delta_{m} \}$. Next, we examine the problem parameters based on the choice of the revision time. In particular, we define the minimum and maximum cost parameters before and after the revision time $t^*$ as follows:
\begin{align*}
\underbar{$a$}^- (t^*) = \min_{n \in {\mathcal{T}}: t_n < t^*} \{a_{n} \}, & \quad \underbar{$a$}^+ (t^*) = \min_{n \in {\mathcal{T}}: t_n \geq t^*} \{a_{n} \}, \\
\bar{a}^- (t^*) = \max_{n \in {\mathcal{T}}: t_n < t^*} \{a_{n} \}, & \quad \bar{a}^+ (t^*) = \max_{n \in {\mathcal{T}}: t_n \geq t^*} \{a_{n} \}.
\end{align*}
For the demand parameters, we let 
\begin{equation*}
\delta^- (t^*) = \max_{m \in {\mathcal{T}}(1, t^* -  1)} \{\delta_{m} \}, \quad  \delta^+(t^*) = \sum_{n \in S_{t^*}} p_n \max_{m \in {\mathcal{T}}(1, t^* -  1) \cup {\mathcal{T}}(n)} \{\delta_{m} \}.
\end{equation*}
Here, the parameter $\delta^- (t^*)$ represents the maximum demand value before the revision time $t^*$. The parameter $\delta^+(t^*)$ corresponds to the expected maximum demand over the tree until the revision time $t^*$ and the subtree rooted at each node of the revision stage. 
To obtain the value of the adaptive two-stage approach, 
we first identify the bounds on the optimum values $v^M$, $v^T$, $v^R (t^*)$.

\begin{prop} \label{PropCapExpBoundAdaptive}
We can derive the following bound for $v^R (t^*)$: 
\begin{equation*}
(\underbar{$a$}^- (t^*) - \underbar{$a$}^+(t^*)) \delta^-(t^*) + \underbar{$a$}^+(t^*) \delta^+(t^*) \leq  v^R (t^*) \leq (\bar{a}^- (t^*) - \bar{a}^+(t^*)) \delta^-(t^*) + \bar{a}^+(t^*) \delta^+(t^*). 
\end{equation*}
\end{prop}

Next, we obtain bounds for the linear programming relaxation of the two-stage equivalent of the formulation \eqref{eq:CapacityExpansionMultiStageReformulation3} using Proposition \ref{PropCapExpBoundAdaptive} by selecting the revision period $t^*$ as 1.

\begin{cor} \label{corCapExpBoundTwoStage}
For $v^T$, we have $a_{*} \delta^* \leq v^T \leq a^* \delta^*$.
\end{cor}  

\begin{prop} \label{PropCapExpBoundMulti} \citep{Huang2009}
For $v^M$, we have $a_{*} \bar \delta \leq v^M \leq a^* \bar \delta$.
\end{prop}

By combining above, we obtain our main result for analyzing the single-resource problem. 

\begin{thm} \label{VatsBoundsTheorem}
We derive the following bounds for  $v^T - v^R (t^*)$ and $v^R (t^*) - v^M$:
\begin{align} \label{eq:subBoundAdjustable}
a_{*} \delta^* - (\bar{a}^- (t^*) - \bar{a}^+(t^*)) \delta^-(t^*) - \bar{a}^+(t^*) \delta^+(t^*) & \leq  v^T - v^R (t^*) \\
& \leq a^* \delta^* - (\underbar{$a$}^- (t^*) - \underbar{$a$}^+(t^*)) \delta^-(t^*) - \underbar{$a$}^+(t^*) \delta^+(t^*), \nonumber
\end{align}
and
\begin{align} \label{eq:subBoundMultistage}
(\underbar{$a$}^- (t^*) - \underbar{$a$}^+(t^*)) \delta^-(t^*) + \underbar{$a$}^+(t^*) \delta^+(t^*) - a^* \bar \delta & \leq  v^R (t^*) - v^M \\
& \leq (\bar{a}^- (t^*) - \bar{a}^+(t^*)) \delta^-(t^*) + \bar{a}^+(t^*) \delta^+(t^*) - a_{*} \bar \delta. \nonumber
\end{align}
\end{thm}

Combining Proposition \ref{PropCapExpBoundAdaptive} and Corollary \ref{corCapExpBoundTwoStage}, we obtain bounds on the VATS for the single-resource problem in \eqref{eq:subBoundAdjustable} ($v^T - v^R (t^*)$). Combining Propositions \ref{PropCapExpBoundAdaptive} and \ref{PropCapExpBoundMulti}, we obtain bounds on the value of the multi-stage model  against adaptive two-stage model for the single-resource problem in \eqref{eq:subBoundMultistage} ($v^R (t^*) - v^M$). These results are critical in deriving insights on the single-resource problem in Section \ref{senstivityVatsSingleResource} and extending them to the capacity expansion planning problem in Section \ref{VATSMultiItem}. 
{\color{black}We note that the lower bounds in Theorem \ref{VatsBoundsTheorem} may potentially take negative values depending on the values of $\{a_n\}_{n \in \mathcal{T}}$, $\{\delta_n\}_{n \in \mathcal{T}}$, and the probabilities associated with the nodes of the scenario tree.}

\subsubsection{Analytical insights on the single-resource problem}
\label{senstivityVatsSingleResource}

In order to gain some insights regarding the effect of the revision point decisions on the performance of the analytical bounds, we consider two cases with respect to the stochasticity of the cost and demand parameters. 

\paragraph{Demand Sensitivity} We first consider the case where the cost parameters $\{a_{n}\}_{n \in {\mathcal{T}}}$ are stationary by not changing over time 
, i.e., {\color{black}$\{a_{n}\}_{n \in {\mathcal{T}}} = a$}. In that case, the bounds \eqref{eq:subBoundAdjustable} and \eqref{eq:subBoundMultistage} reduce to the following expressions: 
{\color{black}
\begin{equation} \label{eq:VATSBounds}
v^T - v^R (t^*) = a (\delta^* - \delta^+(t^*)), \quad v^R (t^*) - v^M = a (\delta^+(t^*) - \bar \delta). 
\end{equation}
}
Consequently, the value of the adaptive formulation depends on how much $\delta^+(t^*)$ differs from the maximum demand value, $\delta^*$. Similarly, the value of the multi-stage stochastic program against the adaptive two-stage approach depends on the difference between $\delta^+(t^*)$ and the maximum average demand, $\bar \delta$. These results highlight the importance of the variability of the demand on the values of the stochastic programs. If there is not much variability in the demand across scenarios, then the three approaches result in similar solutions as the {\color{black}values of the expressions} in \eqref{eq:VATSBounds} go to zero. As variability of the scenarios increases, the corresponding {\color{black}values might increase} accordingly. 

Furthermore, the bounds \eqref{eq:VATSBounds} demonstrate that the values of the adaptive approach might be highly dependent to the selection of $t^*$ value. In particular, best performance gains for the adaptive approach, compared to two-stage and multi-stage cases, can be obtained when $\delta^+(t^*)$ is minimized. By this way, we obtain the solution of the adaptive program which gives the least loss compared to the multi-stage stochastic programs, and  the most gain compared to the two-stage stochastic programs, with respect to the proposed bounds. Let us denote the best revision point in terms of the demand bounds as $t^{DB} := \argmin_{2 \leq t \leq T} \delta^+(t)$. 

We illustrate the effect of the revision times on the bound values \eqref{eq:VATSBounds} for an instance of the single-resource problem described in Appendix \ref{appendixInstanceData}. 
Specifically, we have unit cost for $\{a_{n}\}_{n \in {\mathcal{T}}}$ values, and $\{\delta_{n}\}_{n \in {\mathcal{T}}}$ values are sampled from a probability distribution. 
We present the bound values in Table \ref{BoundComparisonTable}, where $\delta^* = 41$ and $\bar \delta = 34.06$, considering the scenario tree in Figure \ref{IllustrativeInstanceDeltaTree}. We observe that $\delta^+(t^*)$ is minimized at period 3, i.e., $t^{DB} = 3$, maximizing the objective gap of the adaptive two-stage approach with two-stage model and minimizing the corresponding gap with multi-stage model. These results demonstrate that we can analytically determine the revision points for the single-resource problem when the cost parameters are equal to each other. 

\begin{table}[H]
\caption{Bound comparison with respect to the revision time.}
\label{BoundComparisonTable}
\center
\begin{tabular}{cccccc}
\toprule
           &          $t^* = 1$ &          $t^* = 2$ &          $t^* = 3$ &          $t^* = 4$ &          $t^* = 5$ \\
\midrule 
     $v^T - v^R (t^*)$   &       0.00 &       1.50 &       {\bf 5.25} &       2.63 &       0.00 \\
\hline 
       $v^R (t^*) - v^M$ &       6.94 &       5.44 &       {\bf 1.69} &       4.31 &       6.94 \\
\hline
$\delta^+(t^*)$ &      41.00 &      39.50 &      {\bf 35.75} &      38.38 &      41.00 \\
 \bottomrule 
\end{tabular}  
\end{table}

Next, we examine how $\delta^+(t^*)$ changes under different $\delta_{n}$ structures. 
Let $\hat t^D$ be the first stage in which we observe the maximum demand value over the scenario tree, i.e., $\hat t^D = \min \{t \in \{1, \cdots, T\}: \delta_{n} = \delta^*, n \in S_t\}$. 

\begin{prop} \label{demandSensitivityPropGeneric}
Under a general demand structure for $\delta_{n}$ with cost values {\color{black}$\{a_{n}\}_{n \in {\mathcal{T}}} = a$},  if $\hat t^D = 1$, then $t^{DB} \in \{2, \cdots, T\}$. Otherwise, $t^{DB} \in \{2, \cdots, \hat t^D\}$.
\end{prop}

For specific forms of demand patterns, we can further refine Proposition \ref{demandSensitivityPropGeneric}. 

\begin{prop} \label{demandSensitivityPropSpecific}
{\color{black}If the demand values $\{\delta_{n}\}_{n \in {\mathcal{T}}}$ are independently distributed with $N_t$ many realizations for each stage $t_n \in \{2, \cdots, T\}$ with cost values $\{a_{n}\}_{n \in {\mathcal{T}}} = a$, then $t^{DB} = \hat t^D$.}
\end{prop}
Combining above, when demand is independently distributed, the best revision time in terms of this analysis is the period with the maximum demand realization. 
Otherwise, under a general demand distribution, the best revision time is between the second period and the period with the maximum demand realization, unless it is at the first period. 
These results imply that if decisions are revised after observing the maximum demand, then it can be too late, leading to higher costs. 

\paragraph{Cost Sensitivity} Next, we consider the case where the demand parameters $\{\delta_{n}\}_{n \in {\mathcal{T}}}$ are stationary by not changing over time 
, i.e., {\color{black}$\{\delta_{n}\}_{n \in {\mathcal{T}}} = \delta$}. 
Appendix \ref{CostSensitivitySectionAppendix} provides an analysis for this case with different cost patterns. 

\subsubsection{VATS for the capacity expansion planning problem}
\label{VATSMultiItem}

In this section, we extend our results on the single-resource subproblem \eqref{eq:CapacityExpansionAdaptiveReformulation} to the capacity expansion planning problem \eqref{eq:CapacityExpansionMultiStage} under a given revision decision for each resource. We derive analytical bounds on the objective of the adaptive two-stage approach in comparison to {\color{black}the} two-stage and multi-stage stochastic models by leveraging Theorem \ref{VatsBoundsTheorem}. To derive the desired bounds, we utilize the linear programming relaxations of the capacity expansion planning problem. 
\textcolor{black}{We refer to the optimum objective function value of the mixed-integer problem \eqref{eq:CapacityExpansionMultiStage} under the multi-stage, adaptive two-stage with revision decisions at $t^*$, and two-stage models as $V^{MS}$, $V^{ATS}(t^*)$, and $V^{TS}$, respectively. In the remainder of this analysis, once $\delta_n$ value is presented, $\hat \delta_n$ value can be computed for the given condensed tree as discussed in the reformulation from \eqref{eq:CapacityExpansionAdaptiveReformulation} to \eqref{eq:CapacityExpansionAdaptiveReformulationCondensed}.}

\begin{prop} \label{capExpansionBoundProp1}
Let $\{y_n^{TLP}\}_{n \in {\mathcal{T}}}$ be the optimal capacity allocation decisions to the linear programming relaxation of the two-stage version of the stochastic capacity expansion planning problem~\eqref{eq:CapacityExpansionMultiStage}. Then, 
\begin{equation}
V^{TS} - V^{ATS}(t^*) \geq \sum_{i = 1}^I (v_i^T - v_i^R (t_i^*) - \max_{n \in {\hat{\mathcal{T}}}_i (t_i^*)} \{\lceil {\color{black} \hat \delta_{in} \rceil - \hat \delta_{in}}\} a_{i1}),
\end{equation}
where $v_i^T$ and $v_i^R (t_i^*)$ are the optimal objective values of the linear programming relaxations of the models \eqref{eq:CapacityExpansionTwoStageReformulation} and \eqref{eq:CapacityExpansionAdaptiveReformulation} under $\delta_{in} = [A_n y_n^{TLP}]_i$.  
\end{prop} 

\begin{prop} \label{capExpansionBoundProp2}
Let $\{y_n^{MLP}\}_{n \in {\mathcal{T}}}$ be the optimal capacity allocation decisions to the linear programming relaxation of the multi-stage stochastic capacity expansion planning problem~\eqref{eq:CapacityExpansionMultiStage}.
Then, 
\begin{equation}
V^{ATS}(t^*) - V^{MS} \leq \sum_{i = 1}^I (v_i^R (t_i^*) - v_i^M + \max_{n \in {\hat{\mathcal{T}}}_i (t_i^*)} \{\lceil {\color{black} \hat \delta_{in} \rceil - \hat \delta_{in}}\} a_{i1}),
\end{equation}
where $v_i^M$ and $v_i^R (t_i^*)$ are the optimal objective values of the linear programming relaxations of the models \eqref{eq:CapacityExpansionMultiStageReformulation3} and \eqref{eq:CapacityExpansionAdaptiveReformulation} under $\delta_{in} = [A_n y_n^{MLP}]_i$.  
\end{prop} 


Using the bound in Proposition \ref{capExpansionBoundProp2} and the relationship 
\textcolor{black}{$V^{MS} \leq V^{ATS}$}, 
we can evaluate the performance of the adaptive two-stage approach in comparison to its value under a given~$t^*$ vector.  

\begin{cor} \label{boundComparisonAdaptiveTrue}
Let $V^{ATS}$ denote the optimal objective of the adaptive two-stage capacity expansion problem. Then,
\begin{equation}
V^{ATS}(t^*) - V^{ATS} \leq \sum_{i = 1}^I (v_i^R (t_i^*) - v_i^M + \max_{n \in {\hat{\mathcal{T}}}_i (t_i^*)} \{\lceil {\color{black} \hat \delta_{in} \rceil - \hat \delta_{in}}\} a_{i1}),
\end{equation}
where $\delta_{in} = [A_n y_n^{MLP}]_i$.
\end{cor}

These theoretical results are further used in the next section for identifying revisions times to develop solution approaches with approximation guarantees.

\section{Solution Methodology}
\label{approxAlgSection}

The adaptive two-stage formulation imposes computational challenges as the revision times of resources are decision variables. 
\textcolor{black}{In this section, we present a two-phase procedure to solve the adaptive two-stage problems by first obtaining the revision times for each decision requiring limited flexibility, and then solving the resulting problem under these revision decisions.  
In particular, we propose three approximation algorithms to solve the adaptive two-stage equivalent of the capacity expansion planning problem \eqref{eq:CapacityExpansionMultiStage} that leverages this two-phase framework and the bound analyses against the two-stage and multi-stage stochastic programming approaches.} 
Our first two algorithms are based on the bounds derived in Propositions \ref{capExpansionBoundProp1} and \ref{capExpansionBoundProp2} 
by selecting the revision points that provide most gain against two-stage and least loss against multi-stage models. 
In Algorithm~\ref{alg:AdaptiveApproxAlg1}, we identify the revision point of each resource by maximizing the lower bound of the gain in objective in comparison to {\color{black}the} two-stage stochastic model. Similarly, in Algorithm~\ref{alg:AdaptiveApproxAlg2}, we determine the revision points by minimizing the upper bound of the loss in objective in comparison to {\color{black}the} multi-stage stochastic model. 
Furthermore, we propose Algorithm \ref{alg:AdaptiveApproxAlg3} by first solving a relaxation of the adaptive two-stage stochastic program to identify the revision points. In the second-phase of all algorithms, we obtain the capacity expansion and allocation decisions under the resulting revision decisions. 

\begin{algorithm}[h]
\small
\caption{Algorithm Two-stage Relax ({\bf TS-Relax})}
\label{alg:AdaptiveApproxAlg1}
\begin{algorithmic}[1]
\STATE Solve the linear programming relaxation of the two-stage version of the stochastic capacity expansion planning problem \eqref{eq:CapacityExpansionMultiStage} and obtain $\{y_n^{TLP}\}_{n \in {\mathcal{T}}}$.
\STATE Compute $\delta_{in} = [A_n y_n^{TLP}]_i$ for all $i = 1, \cdots, I$, $n \in {\mathcal{T}}$.
\FORALL{$i =1, \dots, I$}
\STATE Find $t_i^*$ that maximizes the lower bound on $v_i^T - v_i^R (t_i^*) - \max_{n \in {\hat{\mathcal{T}}}_i (t_i^*)} \{\lceil  {\color{black}\hat \delta_{in} \rceil - \hat \delta_{in}}\} a_{i1}$ using the bounds~\eqref{eq:subBoundAdjustable}. \label{tStarSolutionStep}
\ENDFOR
\STATE Solve the adaptive two-stage version of the stochastic capacity expansion planning problem \eqref{eq:CapacityExpansionMultiStage} for $\{x_n, y_n\}_{n \in {\mathcal{T}}}$ given the $\{t_i^*\}_{i = 1}^I$ values. 
\end{algorithmic}
\end{algorithm}


\begin{algorithm}[h]
\caption{Algorithm Multi-stage Relax ({\bf MS-Relax})}
\label{alg:AdaptiveApproxAlg2}
\begin{algorithmic}[1]
\small
\STATE Solve the linear programming relaxation of the multi-stage stochastic capacity expansion planning problem \eqref{eq:CapacityExpansionMultiStage} and obtain $\{y_n^{MLP}\}_{n \in {\mathcal{T}}}$.
\STATE Compute $\delta_{in} = [A_n y_n^{MLP}]_i$ for all $i = 1, \cdots, I$, $n \in {\mathcal{T}}$.
\FORALL{$i =1, \dots, I$}
\STATE Find $t_i^*$ that minimizes the upper bound on $v_i^R (t_i^*) - v_i^M + \max_{n \in {\hat{\mathcal{T}}}_i (t_i^*)} \{\lceil  {\color{black} \hat \delta_{in} \rceil - \hat \delta_{in}}\}a_{i1}$ using the bounds \eqref{eq:subBoundMultistage}. \label{tStarSolutionStep2}
\ENDFOR
\STATE Solve the adaptive two-stage version of the stochastic capacity expansion planning problem \eqref{eq:CapacityExpansionMultiStage} for $\{x_n, y_n\}_{n \in {\mathcal{T}}}$ given the $\{t_i^*\}_{i = 1}^I$ values. 
\end{algorithmic}
\end{algorithm}

\begin{algorithm}[h]
\small
\caption{Algorithm Adaptive Two-stage Relax ({\bf ATS-Relax})}
\label{alg:AdaptiveApproxAlg3}
\begin{algorithmic}[1]
\STATE Solve a relaxation of the adaptive two-stage version of the stochastic capacity expansion planning problem~\eqref{eq:CapacityExpansionMultiStage} where $\{x_n\}_{n \in {\mathcal{T}}}$ decisions are continuous \textcolor{black}{and  $\{r_{it}\}$ decisions are binary.  Obtain the revision decisions and name them as $\{r^{ALP}_{it}\}$ for $i = 1, \cdots, I$, $t = 1, \cdots, T$.} 
\STATE Let $t_i^* = \sum_{t = 1}^T t \> r^{ALP}_{it}$ for all $i = 1, \cdots, I$.
\STATE Solve the adaptive two-stage version of the stochastic capacity expansion planning problem \eqref{eq:CapacityExpansionMultiStage} for  $\{x_n, y_n\}_{n \in {\mathcal{T}}}$ given the $\{t_i^*\}_{i = 1}^I$ values. 
\end{algorithmic}
\end{algorithm}

We note that Corollary \ref{boundComparisonAdaptiveTrue} provides an upper bound to compare the objectives of the true adaptive two-stage program with the adaptive program under a given revision point. Using this result, we can demonstrate the approximation guarantees of Algorithms \ref{alg:AdaptiveApproxAlg1} - \ref{alg:AdaptiveApproxAlg3} according to their choices of revision decisions. 
We can further improve this bound for Algorithm~\ref{alg:AdaptiveApproxAlg3} as follows.

\begin{prop} \label{approxGuarenteeProposition}
Let $V^{ATS-Relax}$ and $t^{ATS-Relax}$ denote the objective and the revision vector of the adaptive two-stage program under the solutions found in Algorithm \ref{alg:AdaptiveApproxAlg3}{\color{black}, respectively}. Then,
\begin{equation} \label{eq:ApproxGuarantee}
V^{ATS-Relax} - V^{ATS} \leq \sum_{i = 1}^I (\max_{n \in {\mathcal{T}}_i (t_i^{ATS-Relax})} \{(\lceil  {\color{black}\hat \delta_{in} \rceil - \hat \delta_{in})}\} a_{i1}),
\end{equation}
where $\delta_{in} = [A_n y_n^{LP}]_i$ and $y^{LP}$ is the capacity allocation decisions found in Step 1 of Algorithm \ref{alg:AdaptiveApproxAlg3}.
\end{prop}


We can obtain a simpler upper bound on the optimality gap of the ATS-Relax Algorithm derived in \eqref{eq:ApproxGuarantee}. Specifically, we have $V^{ATS-Relax} - V^{ATS} \leq \sum_{i = 1}^I a_{i1}$ as $\lceil {\color{black} \hat \delta_{in} \rceil - \hat \delta_{in}} \leq 1$ for all resources $i$ and nodes $n$.  This optimality gap is irrespective of the number of stages, scenario tree structure and problem data, but only depends on the capacity acquisition cost of each resource at the first stage. {\color{black}This results in a constant approximation guarantee and asymptotic convergence of the ATS-Relax Algorithm to the true problem as $\lim_{T \to \infty} (V^{ATS-Relax} - V^{ATS})/T \leq \lim_{T \to \infty}  \sum_{i = 1}^I a_{i1}/T =  0$.}

\section{Computational Analysis}
\label{scComputations}

In this section, we provide a comprehensive framework to evaluate the performance of the proposed adaptive two-stage stochastic programming approach over a generation capacity expansion planning problem. 
In Section \ref{experimentalSetupDetails}, we present the experimental setup and the details of the problem formulation. In Section \ref{experimentalResults}, we provide an extensive computational study by demonstrating the value of the proposed approach and illustrating the performances of the proposed solution algorithms on various scenario tree structures. 
We also present a detailed optimal solution for a specific generation expansion problem and provide insights regarding the proposed approach. 

\subsection{Experimental Setup and Optimization Model}
\label{experimentalSetupDetails}

Generation capacity expansion planning is a well-studied problem in {\color{black}the} literature to determine the acquisition and capacity allocation decisions of different types of generation resources over a long-term planning horizon. Our aim is to optimize the yearly investment decisions of different generation resources while producing energy within the available capacity of each generation resource and satisfying the overall system demand in each subperiod. We consider four types of subperiods within a year, namely peak,
shoulder, off-peak and base, depending on their demand level. 
We consider six different types of generation resources for investment decisions, namely nuclear, coal, natural gas-combined cycle (NG-CC), natural gas-gas turbine (NG-GT), wind and solar. We utilize the data set presented in \citet{Min2018}, which is based on a technical report of the U.S. Department of Energy \citep{Black2012}. This data set is used for computing capacity amounts and costs associated with each type of generation resource.

Using the predictions in \citet{Liu2018}, we consider a 10\% reduction in the capacity acquisition costs of the renewable generation resources in each year. Additionally, we assume 15\% and 10\% yearly increases in total for fuel prices and operating costs of natural gas and coal type generation resources, respectively. For other types of generators, we take into account 3\% increase in fuel prices \citep{Min2018}. 
We assume that the generators are available for production in the period that the acquisition decision is made, similar to \citet{Jin2011} and \citet{Zou2018}. 

The source of uncertainty of our model is the demand level of each subperiod 
throughout the planning horizon. We adopt the procedure presented in \citet{Singh2009} for constructing the scenario tree. At the beginning of the planning horizon, we start with an initial demand level for each subperiod, using the values in \citet{Min2018}. Then, we randomly generate a demand increase multiplier for each node,  except the root node, 
to estimate that node's demand value based on its ancestor node's demand. The details of the scenario tree generation algorithm 
are presented in Algorithm \ref{alg:scenarioTreeAlgorithm} in Appendix \ref{scenarioTreeSectionAppendix}. 

Following the studies \citet{Singh2009} and \citet{Zou2018}, we represent the generation capacity expansion planning problem as a stochastic program. We reformulate this problem as an adaptive two-stage program, where the capacity acquisition decisions can be revised once within the planning horizon. \textcolor{black}{Using an analogous formulation to \eqref{eq:CapacityExpansionMultiStage}, we use set $\mathcal{I}$ for the set of generation types for expansion, and set $\mathcal{K}$ for the subperiod types with demand within each period. The parameters and decision variables are defined as follows:}

{\small
\noindent {\bf Parameters:}
\vspace{-3mm}
\begin{flalign*}
& T && \quad \text{Number of periods} & \textcolor{white}{dummydum}\\
& n^0_{i} && \quad \text{Number of generation resource type $i$ at the beginning of planning} \\
& m^{max}_{i} && \quad \text{Maximum capacity of a type $i$ generation resource in MW} \\
& m'^{max}_{i} && \quad \text{Maximum effective capacity of a type $i$ generation resource in MW} \\
& l_i && \quad \text{Peak contribution ratio of a type $i$ generation resource} \\
& c_{in} && \quad \text{Capacity acquisition cost per MW of a generation resource type $i$ at node $n$} \\ 
& f_{in} && \quad \text{Fix operation and maintenance cost per MW of a generation resource type $i$ at node $n$} \\
& g_{in} && \quad \text{Fuel price per MWh of a generation resource type $i$ at node $n$} \\
& k_{in} && \quad \text{Generation cost per MWh of a generation resource type $i$ at node $n$} \\
& h_{kt} && \quad \text{Number of hours of subperiod $k$ in period $t$} \\
& d_{kn} && \quad \text{Hourly demand at subperiod $k$ at node $n$ in MW} \\
& w && \quad \text{Penalty per MWh for demand curtailment} \\
& r && \quad \text{Yearly interest rate} 
\end{flalign*}
\noindent {\bf Decision variables:}
\vspace{-3mm}
\begin{flalign*}
& x_{in} && \quad \text{Number of generation resource type $i$ acquisition at node $n$} & \textcolor{white}{dummydum} \\
& u_{ikn} && \quad \text{Hourly generation amount of generator type $i$ at subperiod $k$ at node $n$ in MW} \\
& v_{kn} && \quad \text{Hourly demand curtailment amount at subperiod $k$ at node $n$ in MW} & \\
& t^*_i && \quad \text{Revision point for acquisition decisions of generation resource type $i$} 
\end{flalign*}
}
We then formulate the generation expansion planning problem as follows:
\begin{subequations} \label{eq:GenerationExpansion}
\begin{alignat}{1}
\min_{x,u,v, t^*} \quad &  \sum_{n \in {\mathcal{T}}} \sum_{i \in \mathcal{I}} p_n \frac{1}{(1+r)^{t_n - 1}} \left(\left(c_{in} + \sum_{t = t_n}^T \frac{f_{in}}{(1+r)^{t - t_n}}\right) m^{max}_{i} x_{in} + \sum_{k \in \mathcal{K}} (g_{in} + k_{in}) h_{k t_n} u_{ikn} \right) \notag \\
& \hspace{8cm} + \sum_{n \in {\mathcal{T}}} \sum_{k \in \mathcal{K}} p_n \frac{1}{(1+r)^{t_n - 1}} w h_{k t_n} v_{kn}\label{eq:GenerationExpansionObj}\\
\text{s.t.} \quad & u_{ikn} \frac{1}{m'^{max}_{i}} \leq n^0_i + \sum_{m \in P(n)} x_{im} \quad \forall n \in {\mathcal{T}}, \quad \forall i \in \mathcal{I}, \quad \forall k \in \mathcal{K}, \label{eq:GenerationExpansionConstr1} \\
& \sum_{i \in I} l_i u_{ikn} + v_{kn} \geq d_{kn} \quad \forall n \in {\mathcal{T}}, \quad \forall k \in \mathcal{K}, \label{eq:GenerationExpansionConstr2}\\
& x_{im} = x_{in} \quad \forall m, n \in S_t, \quad t < t_i^*, \quad \forall i \in \mathcal{I}, \label{eq:GenerationExpansionConstrAdapt1} \\
& x_{im} = x_{in} \quad \forall m, n \in S_t \cap {\mathcal{T}}(j), \quad j \in S_{t_i^*}, \quad t \geq t_i^*, \quad \forall i \in \mathcal{I}, \label{eq:GenerationExpansionConstrAdapt2} \\
& x_{in} \in \mathbb{Z}_{+}, u_{ikn}, v_{kn} \in \mathbb{R}_{+} \quad \forall i \in \mathcal{I}, k \in \mathcal{K}, \forall n \in \mathcal{T}, \quad t^*_i \in \{1, \cdots, T\} \quad \forall i \in \mathcal{I}.
\end{alignat}
\end{subequations}

The objective function \eqref{eq:GenerationExpansionObj} aims to minimize the expected costs of capacity acquisition, allocation and demand curtailment. In order to compute the capacity acquisition costs, we consider the building cost of the generation unit in addition to its maintenance and operating costs for the upcoming periods. For computing the cost of capacity allocation decisions, we use fuel prices and production cost associated with each type of generation resource. All of the costs are then discounted to the beginning of the planning horizon. Constraint \eqref{eq:GenerationExpansionConstr1} ensures that the production amount in each subperiod is restricted by the total available capacity for each type of generation resource. Constraint \eqref{eq:GenerationExpansionConstr2} guarantees that system demand is satisfied. 
\textcolor{black}{As the provided formulation \eqref{eq:GenerationExpansion} is a special form of the capacity expansion problem studied \eqref{eq:CapacityExpansionMultiStage}, constraint \eqref{eq:GenerationExpansionConstr1} corresponds to constraint \eqref{eq:CapacityExpansionMultiStageConstr1}, and constraint \eqref{eq:GenerationExpansionConstr2} refers to constraint \eqref{eq:CapacityExpansionMultiStageConstr2}}. 
Constraints \eqref{eq:GenerationExpansionConstrAdapt1} and \eqref{eq:GenerationExpansionConstrAdapt2} represent the adaptive two-stage relationship for the capacity acquisition decisions such that the acquisition decision of each resource type $i$ can be revised at $t_i^*$, similar to constraints of problem \eqref{eq:CapacityExpansionAdaptiveReformulation}. We note that the presented formulation is nonlinear, however it can be reformulated as a mixed-integer linear program by defining additional variables for the revision decisions, as shown in (\ref{adaptiveTwostageFormulation}). In our computational experiments, we solve the resulting mixed-integer linear program. 
\vspace{-3mm}
\subsection{Computational Results}
\label{experimentalResults}

In this section, 
we first demonstrate the value of the adaptive two-stage approach in comparison to two-stage stochastic programming under different scenario tree structures and demand characteristics. Secondly, we examine the performances of the different algorithms introduced in Section~\ref{approxAlgSection} for solving the adaptive two-stage problem, along with their comparisons to the existing stochastic programming approaches. Finally, we analyze a generation capacity expansion plan under the adaptive two-stage approach to discuss its practical implications. 

To construct our computational testbed, we generate scenario trees for demand values using Algorithm \ref{alg:scenarioTreeAlgorithm} {\color{black}in Appendix \ref{scenarioTreeSectionAppendix}} when number of branches at each period, namely $M$, is equal to 2 or 3. We examine scenario trees with number of stages $T$ ranging from 3 to 8. 
We solve the optimization problem \eqref{eq:GenerationExpansion} under five different randomly generated scenario trees for each $(M,T)$ pair. 
We report the average of these replications to represent the performances of the proposed methods under various instances. \textcolor{black}{We conduct our experiments on an Intel Core i7 2.60GHz machine with 16 GB RAM. We implement the algorithms in Python using Gurobi 9 with a relative optimality gap of 0.1\% and time limit of 2 hours.}

\subsubsection{Value of adaptive two-stage approach}
\label{VATSResults}

To demonstrate the performance of the adaptive two-stage approach in comparison to two-stage stochastic programming, we define the relative value of the adaptive two-stage approach (RVATS) by extending the definition of VATS. Specifically, we let 
\begin{equation} \label{RVATSEquation}
\text{RVATS (\%)} = \frac{(V^{TS} - V^{ATS})}{V^{TS}} \times 100\%, 
\end{equation}
where $V^{TS}$ and $V^{ATS}$ are the objective values corresponding to the two-stage model and adaptive two-stage model, respectively. 
{\color{black}To determine the value of the fully flexible approach, we also introduce the relative value of the multi-stage approach (RVMS) by defining $\text{RVMS (\%)} = (V^{TS} - V^{MS})/V^{TS} \times 100\%$. Furthermore, we define another performance metric that compares the relative value of the adaptive two-stage in comparison to both two-stage and multi-stage stochastic programming as follows:
\begin{equation} \label{RVATSAgainstMSTS}
\frac{(V^{TS} - V^{ATS})}{(V^{TS} - V^{MS})} \times 100\%. \end{equation}
}
{\color{black}First, we illustrate the RVATS under various scenario tree structures in Figure \ref{VATSComparisonGraph}, in comparison to RVMS. Figures \ref{VATSComparisonGraph2Branch} and \ref{VATSComparisonGraph3Branch} show the behavior of the adaptive two-stage and multi-stage approaches under scenario trees with 2 and 3 branches, each with two different demand patterns and averaged over five  instances, respectively. For constructing the demand patterns, we examine the cases with an increasing variance and mean for demand multipliers. Specifically, we consider the demand multipliers in Algorithm \ref{alg:scenarioTreeAlgorithm} as $\underline{\alpha}_t = 1.00$ and $\overline{\alpha}_t= 1.20 + \Gamma t$ for all stages $t = 2, \cdots, 8$ where $\Gamma = 0.05, 0.1$.
To compute RVATS, {\color{black}TS and MS problems are solved to optimality in all instances, and} we utilize $V^{ATS}$ for all {\color{black}2-branch tree instances and 3-branch tree instances with up to 7 stages}. 
For 3-branch cases with 8 stages, we report a lower bound on the RVATS due to the computational difficulty of solving ATS to optimality, as discussed in Section \ref{scAnalyticalAnalysis}. 
In particular, we replace $V^{ATS}$ in Equation \textcolor{black}{\eqref{RVATSEquation}} with $V^{TS-Relax}$, where $V^{TS-Relax}$ represents the objective value of the adaptive two-stage model under the solution of the Algorithm TS-Relax.}

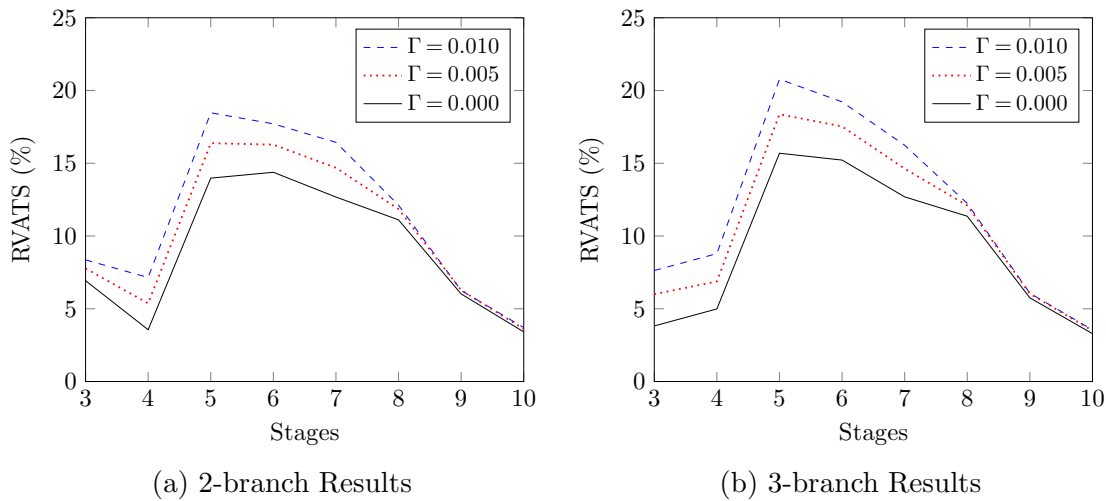
\begin{figure}[h]
\vspace{-3mm}
\caption{Value of Adaptive Two-stage on Instances with Different Variability {\color{black}based on RVATS (\%)}}
\label{VATSComparisonGraph}
\centering 
\begin{subfigure}{.45\textwidth}
\begin{tikzpicture}[scale=0.75]
\begin{axis}[
    xlabel={Stages},
    ylabel={Relative Value (\%)},
    xmin=3, xmax=8,
    ymin=-5, ymax=25,
    xtick={3,4,5,6,7,8},
    ytick={0,5,10,15,20,25}, 
    legend pos=south east,
    grid style=dashed,
]

\addplot[
    color=blue, dashed,
    ]
    table {
    3	11.33396492
4	13.25767889
5	13.91531954
6	15.40805139
7	17.20565538
8	19.13742245
    }; 
\addlegendentry{\small $\Gamma = 0.10$ (RVATS)}

\addplot[
    color=blue,
    ]
    table {
3	11.82
4	14.01
5	14.90
6	16.66
7	18.50
8	20.53
    }; 
\addlegendentry{\small $\Gamma = 0.10$ (RVMS)}

\addplot[
    color=red,thick,dotted,
    ]
    table {
3	9.212704581
4	10.34298908
5	10.08647257
6	10.87614425
7	11.68760136
8	12.5586427
    }; 
\addlegendentry{\small $\Gamma = 0.05$ (RVATS)}

\addplot[
    color=red,thick,
    ]
    table {
3	9.834499692
4	11.1043859
5	11.10582001
6	12.08150453
7	13.00968228
8	13.92512865
    }; 
\addlegendentry{\small $\Gamma = 0.05$ (RVMS)}

\end{axis}
\end{tikzpicture}
\caption{\small 2-branch Results}
\label{VATSComparisonGraph2Branch}
\end{subfigure}
\begin{subfigure}{.45\textwidth}
\begin{tikzpicture}[scale=0.75]
\begin{axis}[
     xlabel={Stages},
    ylabel={Relative Value (\%)},
    xmin=3, xmax=8,
    ymin=-5, ymax=25,
    xtick={3,4,5,6,7,8},
    ytick={0,5,10,15,20,25}, 
    legend pos=south east,
    grid style=dashed,
]

\addplot[
    color=blue, dashed,
    ]
    table {
3	14.70027562
4	14.10541264
5	14.8823642
6	16.33399973
7	16.99927043
8	15.56198641
    }; 
\addlegendentry{\small $\Gamma = 0.10$ (RVATS)}

\addplot[
    color=blue,
    ]
    table {
3	15.52
4	14.97
5	16.01
6	17.68
7	18.50
8	20.38
    }; 
\addlegendentry{\small $\Gamma = 0.10$ (RVMS)}

\addplot[
    color=red,thick,dotted,
    ]
    table {
3	12.08568672
4	11.2796992
5	10.93155976
6	11.59567586
7	11.6014021
8	10.36679351
    }; 
\addlegendentry{\small $\Gamma = 0.05$ (RVATS)}

\addplot[
    color=red,thick,
    ]
    table {
3	12.97007823
4	12.17969338
5	12.11234091
6	12.99300019
7	13.11734038
8	13.9104005
    }; 
\addlegendentry{\small $\Gamma = 0.05$ (RVMS)}
\end{axis}
\end{tikzpicture}
\caption{\small 3-branch Results}
\label{VATSComparisonGraph3Branch}
\end{subfigure}
\end{figure}

{\color{black}We observe RVATS between 9-20\% over 2-branch and 3-branch scenario trees, demonstrating the significant gain of revising decisions over two-stage stochastic programming policies. 
Furthermore, we present the relative value of the adaptive two-stage approach by using the performance metric \eqref{RVATSAgainstMSTS} in Table \ref{tab:RelativeValueComparisonInstances} over each  instance separately.
Although the adaptive two-stage framework enables one revision time for each resource, it obtains a similar relative value in comparison to the multi-stage setting, with performance metrics greater than 85\% over majority of the instances. 
Thus, despite this limited flexibility, the adaptive two-stage approach is consistently much closer to the multi-stage approach than the two-stage, demonstrating the value of optimizing the revision decisions that can capture best time to revise actions, instead of adopting the two-stage approach. 
\vspace{-5mm}
\begin{table}[h]
\centering
\small
\begin{tabular}{cccrrrrrr}
\toprule
           &            &            & \multicolumn{6}{c}{{\color{black}number of stages}} \\
\cline{4-9}
\# of branches &     $\Gamma$ & instance no &          3 &          4 &          5 &          6 &          7 &          8           \\
\midrule 
         2 &       0.05 &          1 &      91.89 &      93.43 &      89.10 &      85.34 &      87.73 &      87.99 \\

           &            &          2 &      91.61 &      93.17 &      89.05 &      90.16 &      89.50 &      90.19 \\

           &            &          3 &      99.77 &      94.75 &      92.58 &      89.78 &      88.64 &      89.41 \\

           &            &          4 &      96.68 &      93.91 &      93.06 &      92.08 &      92.20 &      92.61 \\

           &            &          5 &      89.89 &      89.96 &      90.12 &      92.43 &      90.80 &      90.32 \\
\hline 
         2 &        0.1 &          1 &      95.61 &      95.14 &      91.26 &      88.97 &      91.19 &      91.21 \\

           &            &          2 &      92.35 &      94.64 &      93.03 &      91.17 &      92.26 &      92.94 \\

           &            &          3 &      98.39 &      95.54 &      93.76 &      92.27 &      91.81 &      92.67 \\

           &            &          4 &      97.99 &      94.71 &      95.44 &      94.18 &      94.99 &      95.16 \\

           &            &          5 &      94.89 &      92.39 &      93.36 &      94.24 &      94.08 &      93.47 \\
\hline 
         3 &       0.05 &          1 &      95.30 &      92.83 &      91.53 &      89.38 &      88.16 &      60.80 \\

           &            &          2 &      93.17 &      92.86 &      90.88 &      90.73 &      89.46 &      64.06 \\

           &            &          3 &      93.65 &      92.84 &      88.38 &      87.90 &      87.38 &      53.25 \\

           &            &          4 &      92.73 &      93.23 &      91.38 &      89.59 &      88.98 &      71.66 \\

           &            &          5 &      90.88 &      91.20 &      89.33 &      88.62 &      88.30 &      57.89 \\
\hline 
         3 &        0.1 &          1 &      95.45 &      94.53 &      93.76 &      92.16 &      91.46 &      62.43 \\

           &            &          2 &      96.17 &      94.31 &      93.72 &      93.85 &      92.85 &      88.35 \\

           &            &          3 &      95.28 &      94.28 &      91.53 &      91.29 &      90.98 &      75.62 \\

           &            &          4 &      94.58 &      94.67 &      94.20 &      92.91 &      92.42 &      73.31 \\

           &            &          5 &      91.99 &      93.27 &      91.44 &      91.64 &      91.74 &      62.92 \\
           \bottomrule 
\end{tabular}  
    \caption{The Value of the Adaptive Two-stage Approach in Comparison to Multi-stage and Two-stage Stochastic Programming using the Performance Metric \eqref{RVATSAgainstMSTS}.}
    \label{tab:RelativeValueComparisonInstances}
    \vspace{-3mm}
\end{table}

We note that the scenario tree has more variability in each stage when the parameter $\Gamma$ increases. We observe higher gain of the adaptive two-stage approach in terms of RVATS, for the scenario trees with larger variability levels. The adaptive two-stage approach can preserve its value against the multi-stage approach for scenario trees with 2-branches even if the number of stages increases. 
When scenario tree has more branches, then fully flexible approaches can obtain more relative value against the two-stage and adaptive two-stage approaches for problems with larger number of stages. This situation arises as allowing one revision results in a more condensed scenario tree, i.e., a more rigid decision structure for scenario trees with higher number of branches. 
}


\subsubsection{Performance of solution algorithms}
\label{algorithmResults}

In this section, we examine the performance of the solution algorithms for the capacity expansion planning problem by comparing them with respect to the existing stochastic programming approaches, along with their computational performances. Tables \ref{2BranchAlgResults} and \ref{3BranchAlgResults} provide the results of the solution algorithms TS-Relax, MS-Relax, and ATS-Relax for 2-branch and 3-branch scenario trees, respectively. \textcolor{black}{Additionally, the abbreviation ``ATS'' stands for the true adaptive two-stage approach when revision times are also decision variables. ``MS'' and ``TS'' represent the multi-stage and two-stage approaches, respectively{\color{black}, which are solved to optimality under their extensive formulations}.
To compare the proposed approach and the corresponding solution algorithms with the existing approaches allowing different levels of flexibility, we report the relative percentage gain of the studied method in comparison to {\color{black}the} two-stage stochastic model's objective under the column ``Gain'', and the relative percentage loss in comparison to multi-stage stochastic model's objective under the column ``Loss''. We note that the relative percentage gain for the ATS is equivalent to the relative value of this mehod, namely RVATS introduced in \eqref{RVATSEquation}. We further clarify that although we report comparisons against the multi-stage approach, it does not provide a feasible solution for the adaptive two-stage models as it presents a fully adaptive solution by neglecting the limited flexibility of the underlying problem setting.} %

{\color{black}We conduct the computational experiments under the demand multipliers  $\underline{\alpha}_t = 1.00$ and $\overline{\alpha}_t= 1.20 + \Gamma t$ for all stages $t = 2, \cdots, 8$ where $\Gamma = 0.1$, and present the average performance over five replications. 
Relative gain and loss percentages are not reported for some larger scale instances with {\color{black}the} ATS approach as they cannot be solved to the optimality within the given time limit.} 
We highlight that the approximation algorithms MS-Relax and ATS-Relax provide optimality gaps for the adaptive two-stage program. In particular, these algorithms obtain a lower bound to the adaptive two-stage program by solving relaxations of the multi-stage model and adaptive two-stage model in Step 1 of the Algorithm \ref{alg:AdaptiveApproxAlg1} and Algorithm \ref{alg:AdaptiveApproxAlg3}. Since both algorithms construct a feasible solution, they provide an upper bound for the desired problem. Combining the lower and upper bounds, we construct optimality gap for the adaptive two-stage approach. The optimality gap values are reported in the column ``\%Gap'' in Tables \ref{2BranchAlgResults} and \ref{3BranchAlgResults}. These values provide performance metrics for the solution algorithms in approximating the adaptive two-stage model. 


\begin{table}[h]
\color{black}{
\centering
\caption{Performance of solution algorithms on 2-branch results.}
\label{2BranchAlgResults} 
\begin{tabular}{r|r|rr|rr|rrr|rrr|r}
\toprule
           &         \multicolumn{1}{|c|}{MS} & \multicolumn{ 2}{|c|}{ATS} & \multicolumn{ 2}{|c|}{TS-Relax} &        \multicolumn{ 3}{|c|}{MS-Relax} &       \multicolumn{ 3}{|c|}{ATS-Relax} &         \multicolumn{1}{|c}{TS} \\
\midrule
     Stage &       \%Gain &      \%Gain  &       \%Loss &      \%Gain  &       \%Loss &      \%Gain  &       \%Loss &        \%Gap &      \%Gain  &       \%Loss &        \%Gap &       \%Loss \\
\hline

         3 &      11.82 &      11.33 &       0.56 &      11.33 &       0.56 &      10.89 &       1.05 &       1.40 &      11.33 &       0.56 &       0.31 &      13.46 \\

         4 &      14.01 &      13.26 &       0.88 &      12.95 &       1.24 &      13.07 &       1.10 &       1.32 &      13.25 &       0.89 &       0.22 &      16.44 \\

         5 &      14.90 &      13.92 &       1.16 &      12.42 &       2.94 &      12.80 &       2.46 &       2.62 &      13.90 &       1.18 &       0.14 &      17.54 \\

         6 &      16.66 &      15.41 &       1.51 &      12.42 &       5.09 &      14.81 &       2.23 &       2.36 &      15.39 &       1.52 &       0.11 &      20.05 \\

         7 &      18.50 &      17.21 &       1.59 &      14.67 &       4.70 &      15.91 &       3.18 &       3.29 &      17.17 &       1.64 &       0.10 &      22.81 \\

         8 &      20.53 &      19.14 &       1.76 &      16.69 &       4.83 &      16.28 &       5.37 &       5.47 &      19.10 &       1.81 &       0.10 &      25.99 \\


\bottomrule
\end{tabular}  
}
\end{table}
\begin{table}[h]
\color{black}{
\centering
\caption{Performance of solution algorithms on 3-branch results.}
\label{3BranchAlgResults} 
\begin{tabular}{r|r|rr|rr|rrr|rrr|r}
\toprule
           &         \multicolumn{1}{|c|}{MS} & \multicolumn{ 2}{|c|}{ATS} & \multicolumn{ 2}{|c|}{TS-Relax} &        \multicolumn{ 3}{|c|}{MS-Relax} &       \multicolumn{ 3}{|c|}{ATS-Relax} &         \multicolumn{1}{|c}{TS} \\
\midrule
     Stage &       \%Gain &      \%Gain  &       \%Loss &      \%Gain  &       \%Loss &      \%Gain  &       \%Loss &        \%Gap &      \%Gain  &       \%Loss &        \%Gap &       \%Loss \\
\hline

         3 &      15.52 &      14.70 &       0.97 &      14.70 &       0.97 &      14.70 &       0.97 &       1.35 &      14.70 &       0.97 &       0.33 &      18.44 \\

         4 &      14.97 &      14.11 &       1.02 &      12.68 &       2.74 &      13.58 &       1.65 &       1.87 &      14.11 &       1.02 &       0.17 &      17.65 \\

         5 &      16.01 &      14.88 &       1.35 &      11.32 &       5.58 &      13.82 &       2.59 &       2.73 &      14.87 &       1.36 &       0.14 &      19.10 \\

         6 &      17.68 &      16.33 &       1.64 &      12.64 &       6.12 &      15.21 &       3.00 &       3.11 &      16.33 &       1.64 &       0.10 &      21.52 \\

         7 &      18.50 &      17.00 &       1.84 &      13.60 &       6.01 &      13.60 &       6.01 &       6.12 &      16.98 &       1.87 &       0.08 &      22.72 \\

         8 &      20.38 &      - &       - &      15.56 &       6.05 &      15.67 &       5.91 &       6.00 &      18.51 &       2.36 &       0.10 &      25.61 \\


\bottomrule
\end{tabular}
}
\end{table}

 \textcolor{black}{The MS approach provides the most gain against the two-stage approach, as it provides a relaxation to the studied setting. 
 ATS approach provides gains similar to the MS approach and its loss against multi-stage approach is between 0.5\% and 2\% over all instances. Although ATS model considers problem settings with partial flexibility, it is powerful in representing the underlying uncertainty by optimizing the best revision time for each resource and obtaining a solution with highest gain and least loss under these limitations.}

\textcolor{black}{Among the solution  algorithms, ATS-Relax has the highest gain and least loss despite of its computational disadvantage. The gain of Algorithm ATS-Relax is very close to the gain of ATS demonstrating the success of the proposed algorithm in approximating the adaptive two-stage approach.} Our computational results also highlight the asymptotic convergence of Algorithm ATS-Relax. 
Specifically, as number of stages increases, the optimality gap provided by Algorithm ATS-Relax becomes closer to zero. \textcolor{black}{Additionally, we note that Algorithms TS-Relax and MS-Relax are not superior to each other as in some settings TS-Relax outperforms MS-Relax, and in some other settings MS-Relax performs better than TS-Relax with respect to the relative gain and loss values.} 

\textcolor{black}{The size of the scenario tree and consequently the problem size increase as the number of stages and the number of branches become larger. We examine the computational performance of the ATS approach along with its approximation algorithms in Figure~\ref{computationalPerformanceFigure}. We observe that the computational complexity of the ATS and ATS-Relax approaches are more sensitive to the instance size in comparison to Algorithms TS-Relax and MS-Relax. 
Although ATS-Relax approach provides asymptomatic convergence guarantee, we observe significant computational benefits of adopting 
Algorithms TS-Relax and MS-Relax as they provide notable speedups in comparison to ATS.} 

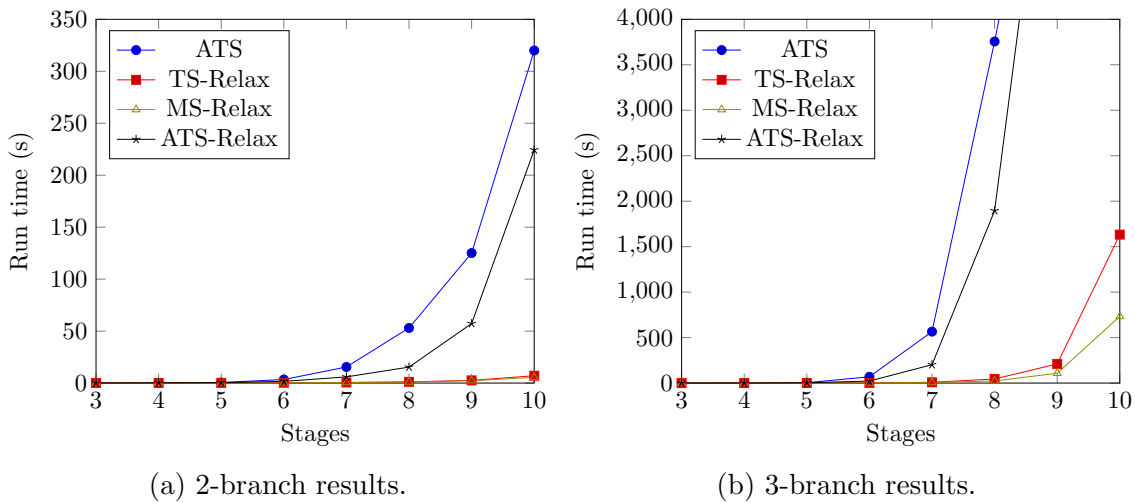
\begin{figure}[h]
\centering
\captionsetup{justification=centering}
\caption{Computational performance of different solution methodologies.}
\label{computationalPerformanceFigure}
\begin{subfigure}[b]{.45\textwidth}
  \centering
\begin{tikzpicture}[scale=0.75]
\begin{axis}[
	xlabel={Stages},
	ylabel={Run time (s)},
	ymin = 0, ymax= 1000,
	xmin = 3, xmax = 8,
	xtick = {3,4,5,6,7,8},
	ytick = {0,200,400,600,800,1000}, 
     legend pos=north west,
]
\addplot coordinates {
(	3	,	0.162712479	)
(	4	,	2.132268953	)
(	5	,	7.236488533	)
(	6	,	23.81427965	)
(	7	,	642.5294951	)
(	8	,	843.1850868)
};

\addplot coordinates {
(	3	,	0.043606949	)
(	4	,	0.202011967	)
(	5	,	0.470918083	)
(	6	,	0.430476952	)
(	7	,	1.312051392	)
(	8	,	3.236927176	)
};

\addplot[mark=triangle,color=olive] coordinates {
(	3	,	0.043569279	)
(	4	,	0.164913321	)
(	5	,	0.269890594	)
(	6	,	0.610247135	)
(	7	,	1.455201244	)
(	8	,	2.617832851	)
};

\addplot coordinates {
(	3	,	0.071635866	)
(	4	,	0.318447399	)
(	5	,	0.624533796	)
(	6	,	2.073785591	)
(	7	,	5.316349554	)
(	8	,	20.40688095	)
};

\legend{ATS, TS-Relax, MS-Relax, ATS-Relax}
\end{axis}
\end{tikzpicture}
  \caption{\small 2-branch results.}
\end{subfigure}
\begin{subfigure}[b]{.45\textwidth}
  \centering
  \begin{tikzpicture}[scale=0.75]
\begin{axis}[
	xlabel={Stages},
	ylabel={Run time (s)},
	ymin = 0, ymax= 9000,
	xmin = 3, xmax = 8,
	xtick = {3,4,5,6,7,8},
	ytick = {0,1000,2000,3000,4000,5000,6000,7000,8000,9000},
     legend pos=north west,
]
\addplot coordinates {
(	3	,	1.917116976	)
(	4	,	11.68902588	)
(	5	,	42.43221598	)
(	6	,	403.7008141	)
(	7	,	6581.678332	)
(	8	,	8778.914289	)  
};

\addplot coordinates {
(	3	,	14.69693469	)
(	4	,	12.68259676	)
(	5	,	11.32350535	)
(	6	,	12.64360915	)
(	7	,	13.60317991	)
(	8	,	15.56198641	)
};

\addplot[mark=triangle,color=olive] coordinates {
(	3	,	0.115408278	)
(	4	,	0.446314907	)
(	5	,	2.218569469	)
(	6	,	4.454554272	)
(	7	,	16.11203346	)
(	8	,	183.5797051	)
};

\addplot coordinates {
(	3	,	0.163484478	)
(	4	,	1.562673664	)
(	5	,	5.944309664	)
(	6	,	29.21268497	)
(	7	,	385.7363405	)
(	8	,	7944.39356 ) 
};

\legend{ATS, TS-Relax, MS-Relax, ATS-Relax}
\end{axis}
\end{tikzpicture}
   \caption{\small 3-branch results.}
\end{subfigure}
\end{figure}

\subsubsection{Discussion on capacity expansion plans}

This section examines a particular instance to analyze the generation capacity expansion plan of the adaptive two-stage model and compare it with that of the two-stage model. \textcolor{black}{Figure \ref{AdaptiveGenerationPlanFigure} illustrates an expansion plan under the adaptive two-stage model over a five-year planning horizon.} For each node of the tree, we show the acquired effective capacity amount of each generation resource type $i \in \{1, \cdots, 6\}$ in the order of nuclear, coal, NG-CC, NG-GT, wind and solar. \textcolor{black}{The revision times of each resource type are selected by the model as 4, 3, 2, 2, 5, 5, 
respectively, and the expansion amounts of those periods are denoted in bold. We report the capacity expansion decisions of a resource type $i$ in node $n$ if period $t_n = t_i^*$ or $x_{in} > 0$. For notational simplicity, we omit the node index in the illustration. 
Figure \ref{TwoStageGenerationPlanFigure} shows the generation expansion plan of the same instance under two-stage model.} 
The initial effective capacities of the resources are 4860, 2618, 2004, 1883, 134, 131 MW, respectively. 

\begin{figure}[h]
\caption{Generation expansion plan \textcolor{black}{with respect to capacity acquisition} under the adaptive two-stage model.}
\label{AdaptiveGenerationPlanFigure}
\centering
\scriptsize
\begin{tikzpicture}[shorten >=1pt,->,draw=black!50, node distance=\layersep, thick, scale=0.3]
    \tikzstyle{every pin edge}=[<-,shorten <=1pt]
    \tikzstyle{neuron}=[minimum size=15pt,inner sep=2pt, draw=black] 
    \tikzstyle{input neuron}=[neuron];
    \tikzstyle{output neuron}=[neuron, fill=white!00, draw=white];
    \tikzstyle{hidden neuron}=[neuron, fill=blue!50];
    \tikzstyle{annot} = [text width=4em, text centered]
    
	\node[input neuron] (I0) at (0,0) {$x_{4}=753$};
	\node[input neuron] (I1) at (5,4) {$\begin{bmatrix} x_{1} \\ x_{2} \\ \mathbf{x_{3}} \\ \mathbf{x_{4}} \end{bmatrix} = \begin{bmatrix} 972 \\ 524 \\ \mathbf{1002} \\ \mathbf{1130} \end{bmatrix}$};
	
	\node[input neuron] (I2) at (5,-4) {$\begin{bmatrix} x_{1} \\ x_{2} \\ \mathbf{x_{3}} \\ \mathbf{x_{4}} \end{bmatrix} = \begin{bmatrix} 972 \\ 524 \\ \mathbf{0} \\ \mathbf{0} \end{bmatrix}$};
	\node[input neuron] (I3) at (13,8.5) {$\begin{bmatrix} x_{1} \\ \mathbf{x_{2}} \\ x_{3} \end{bmatrix} = \begin{bmatrix} 1944 \\	\mathbf{2094} \\ 1002
 \end{bmatrix}$};
	\node[input neuron] (I4) at (13,2.5) {$\begin{bmatrix} x_{1} \\ \mathbf{x_{2}} \\ x_{3} \end{bmatrix} = \begin{bmatrix} 1944 \\	\mathbf{0} \\ 1002
 \end{bmatrix}$};
	\node[input neuron] (I5) at (13,-2.5) {$\begin{bmatrix} x_{1} \\ \mathbf{x_{2}} \\ x_{3} \end{bmatrix} = \begin{bmatrix} 1944 \\	\mathbf{0} \\ 1002
 \end{bmatrix}$};
	\node[input neuron] (I6) at (13,-8.5) {$\begin{bmatrix} x_{1} \\ \mathbf{x_{2}} \\ x_{3} \end{bmatrix} = \begin{bmatrix} 1944 \\	\mathbf{0} \\ 1002
 \end{bmatrix}$};
	\node[input neuron] (I7) at (21,13.5) {$\mathbf{x_{1}}=\mathbf{1944}$};
	\node[input neuron] (I8) at (21,11.5) {$\mathbf{x_{1}}=\mathbf{1944}$};
	\node[input neuron] (I9) at (21,7.5) {$\begin{bmatrix} \mathbf{x_{1}} \\ x_{2} \end{bmatrix} = \begin{bmatrix} \mathbf{1944} \\	2094 \end{bmatrix}$};
	\node[input neuron] (I10) at (21,3.5) {$\begin{bmatrix} \mathbf{x_{1}} \\ x_{2} \end{bmatrix} = \begin{bmatrix} \mathbf{1944} \\	2094 \end{bmatrix}$};
	\node[input neuron] (I11) at (21,-0.5) {$\begin{bmatrix} \mathbf{x_{1}} \\ x_{2} \\ x_{3} \\ x_{4} \end{bmatrix} = \begin{bmatrix} \mathbf{1944} \\	2094 \\ 1002 \\ 1130
 \end{bmatrix}$};
	\node[input neuron] (I12) at (21,-5.5) {$\begin{bmatrix} \mathbf{x_{1}} \\ x_{2} \\ x_{3} \\ x_{4} \end{bmatrix} = \begin{bmatrix} \mathbf{972} \\	2094 \\ 1002 \\ 1130
 \end{bmatrix}$};
	\node[input neuron] (I13) at (21,-11.5) {$\begin{bmatrix} \mathbf{x_{1}} \\ x_{2} \\ x_{3} \\ x_{4} \end{bmatrix} = \begin{bmatrix} \mathbf{1944} \\	524 \\	1002	\\ 1130
 \end{bmatrix}$};
	\node[input neuron] (I14) at (21,-16.5) {$\begin{bmatrix} \mathbf{x_{1}} \\ x_{2} \\ x_{3} \\ x_{4} \end{bmatrix} = \begin{bmatrix} \mathbf{0} \\	524 \\	1002	\\ 1130
 \end{bmatrix}$};
	\node[input neuron] (I15) at (30,19.25) {$\begin{bmatrix} \mathbf{x_{5}} \\ \mathbf{x_{6}}  \end{bmatrix} = \begin{bmatrix} \mathbf{3000} \\	\mathbf{4800} \end{bmatrix}$};
	\node[input neuron] (I16) at (30,16.25) {$\begin{bmatrix} \mathbf{x_{5}} \\ \mathbf{x_{6}} \end{bmatrix} = \begin{bmatrix} \mathbf{1929} \\	\mathbf{4800} \end{bmatrix}$};
	\node[input neuron] (I17) at (30,13.25) {$\begin{bmatrix} \mathbf{x_{5}} \\ \mathbf{x_{6}} \end{bmatrix} = \begin{bmatrix} \mathbf{1356} \\	\mathbf{4800} \end{bmatrix}$}; 
	\node[input neuron] (I18) at (30,10.25) {$\begin{bmatrix} \mathbf{x_{5}} \\ \mathbf{x_{6}} \end{bmatrix} = \begin{bmatrix} \mathbf{0} \\	\mathbf{0} \end{bmatrix}$};
	\node[input neuron] (I19) at (30,7.25) {$\begin{bmatrix} \mathbf{x_{5}} \\ \mathbf{x_{6}} \end{bmatrix} = \begin{bmatrix} \mathbf{2568} \\	\mathbf{4800} \end{bmatrix}$};
	\node[input neuron] (I20) at (30,4.25) {$\begin{bmatrix} \mathbf{x_{5}} \\ \mathbf{x_{6}} \end{bmatrix} = \begin{bmatrix} \mathbf{0} \\	\mathbf{0} \end{bmatrix}$};
	\node[input neuron] (I21) at (30,1.25) {$\begin{bmatrix} \mathbf{x_{5}} \\ \mathbf{x_{6}} \end{bmatrix} = \begin{bmatrix} \mathbf{0} \\	\mathbf{0} \end{bmatrix}$};
	\node[input neuron] (I22) at (30,-1.75) {$\begin{bmatrix} \mathbf{x_{5}} \\ \mathbf{x_{6}} \end{bmatrix} = \begin{bmatrix} \mathbf{0} \\	\mathbf{0} \end{bmatrix}$};
	\node[input neuron] (I23) at (30,-4.75) {$\begin{bmatrix} \mathbf{x_{5}} \\ \mathbf{x_{6}} \end{bmatrix} = \begin{bmatrix} \mathbf{3000} \\	\mathbf{4800} \end{bmatrix}$};
	\node[input neuron] (I24) at (30,-7.75) {$\begin{bmatrix} \mathbf{x_{5}} \\ \mathbf{x_{6}} \end{bmatrix} = \begin{bmatrix} \mathbf{0} \\	\mathbf{0} \end{bmatrix}$};
	\node[input neuron] (I25) at (30,-10.75) {$\begin{bmatrix} \mathbf{x_{5}} \\ \mathbf{x_{6}}  \end{bmatrix} = \begin{bmatrix} \mathbf{0} \\	\mathbf{0} \end{bmatrix}$};
	\node[input neuron] (I26) at (30,-13.75) {$\begin{bmatrix} \mathbf{x_{5}} \\ \mathbf{x_{6}} \end{bmatrix} = \begin{bmatrix} \mathbf{0} \\	\mathbf{0} \end{bmatrix}$};
	\node[input neuron] (I27) at (30,-17) {$\begin{bmatrix} x_{2} \\ \mathbf{x_{5}} \\ \mathbf{x_{6}}  \end{bmatrix} = \begin{bmatrix} 1571 \\\mathbf{0} \\ \mathbf{427} \end{bmatrix}$};
	\node[input neuron] (I28) at (30,-21) {$\begin{bmatrix} x_{2} \\ \mathbf{x_{5}} \\ \mathbf{x_{6}}  \end{bmatrix} = \begin{bmatrix} 1571 \\ \mathbf{0} \\	\mathbf{0} \end{bmatrix}$};
	\node[input neuron] (I29) at (30,-25) {$\begin{bmatrix} x_{2} \\ \mathbf{x_{5}} \\ \mathbf{x_{6}}  \end{bmatrix} = \begin{bmatrix} 1571 \\ \mathbf{0} \\	\mathbf{0} \end{bmatrix}$};
	\node[input neuron] (I30) at (30,-29) {$\begin{bmatrix} x_{2} \\ \mathbf{x_{5}} \\ \mathbf{x_{6}}  \end{bmatrix} = \begin{bmatrix} 1571 \\ \mathbf{0} \\	\mathbf{0} \end{bmatrix}$};
	
            \path (I0) edge (I1);
            \path (I0) edge (I2);
            \path (I1) edge (I3);
            \path (I1) edge (I4);
            \path (I2) edge (I5);
            \path (I2) edge (I6);
            \path (I3) edge (I7);
            \path (I3) edge (I8);
            \path (I4) edge (I9);
            \path (I4) edge (I10);
            \path (I5) edge (I11);
            \path (I5) edge (I12);
            \path (I6) edge (I13);
            \path (I6) edge (I14);
            \path (I7) edge (I15);
            \path (I7) edge (I16);
            \path (I8) edge (I17);
            \path (I8) edge (I18);
            \path (I9) edge (I19);
            \path (I9) edge (I20);
            \path (I10) edge (I21);
            \path (I10) edge (I22);
            \path (I11) edge (I23);
            \path (I11) edge (I24);
            \path (I12) edge (I25);
            \path (I12) edge (I26);
            \path (I13) edge (I27);
            \path (I13) edge (I28);
            \path (I14) edge (I29);
            \path (I14) edge (I30);
    
\end{tikzpicture}
\end{figure}
\vspace{-1mm}
\begin{figure}[h]
\caption{Generation expansion plan \textcolor{black}{with respect to capacity acquisition} under the two-stage model.}
\label{TwoStageGenerationPlanFigure}
\centering
\scriptsize
\begin{tikzpicture}[shorten >=1pt,->,draw=black!50, node distance=\layersep, thick, scale=0.3]
    \tikzstyle{every pin edge}=[<-,shorten <=1pt]
    \tikzstyle{neuron}=[minimum size=15pt,inner sep=2pt, draw=black] 
    \tikzstyle{input neuron}=[neuron];
    \tikzstyle{output neuron}=[neuron, fill=white!00, draw=white];
    \tikzstyle{hidden neuron}=[neuron, fill=blue!50];
    \tikzstyle{annot} = [text width=4em, text centered]  
	
	\node[input neuron] (I0) at (0,0) {$x_4 = 753$};
	\node[input neuron] (I1) at (7,0) {$\begin{bmatrix} x_{1} \\ x_{2} \\ x_{3} \\ x_{4} \end{bmatrix} = \begin{bmatrix} 972 \\ 524 \\ 501 \\ 1130 \end{bmatrix}$};
	\node[input neuron] (I2) at (16,0) {$\begin{bmatrix} x_{1} \\ x_{2} \\ x_{3} \end{bmatrix} = \begin{bmatrix} 2916 \\ 1047 \\ 1503 \end{bmatrix}$};
	\node[input neuron] (I3) at (24.5,0) {$\begin{bmatrix} x_{1} \\ x_{2} \end{bmatrix} = \begin{bmatrix} 972 \\ 1047 \end{bmatrix}$};
	\node[input neuron] (I4) at (31,0) {$x_6 = 4800$};
		
            \path (I0) edge (I1);
            \path (I1) edge (I2);
            \path (I2) edge (I3);
            \path (I3) edge (I4);
    
\end{tikzpicture}
\end{figure}
\textcolor{black}{We observe that the adaptive two-stage approach provides more flexibility compared to the two-stage by allowing an update for each generation type's expansion plan. Since the  two-stage approach is less adaptive to the uncertainty, it brings forward the expansion times of the resources to satisfy the overall demand of each stage. We note that the expected total capacity expansion and allocation cost of the adaptive two-stage model 
results in 12.30\% relative gain of this approach against its two-stage counterpart. Although the multi-stage approach does not provide a feasible solution due to having a fully flexible investment plan for this setting requiring limited flexibility, we note that the relative gain of the multi-stage approach is 12.87\%, which is very close to the relative gain of the adaptive approach.} 

\textcolor{black}{As the investment costs of the renewable resources decrease over time, both of the models tend to expand wind and solar generation capacities at later times in the planning. Due to this decreasing  investment costs and increasing variability of demand in the scenario tree, the revision time of wind and solar capacities are in $5^{th}$ periods to adjust the expansion decisions based on the underlying demand. Consequently, these renewable resources' capacity acquisition amounts vary significantly in the later periods depending on the realized uncertainty, resulting in significant differences with the results of the two-stage model. 
The fuel prices and operating costs associated with traditional generation types increase over time. Thus, the capacities of these types of resources are mostly expanded in the earlier periods of the planning horizon. To adapt to the demand uncertainty, we observe revisions of expansion decisions at $2^{nd}$ periods for natural gas type generation resources, and $4^{th}$ and $3^{rd}$ periods for nuclear and coal based energy generation.}

\section{Conclusion}
\label{scConclusion}

In this paper, we propose a stochastic optimization {\color{black}framework} that determines the best time to revise decisions for problems allowing partial flexibility. 
\textcolor{black}{We first present a generic formulation for the proposed adaptive two-stage approach, 
and highlight the importance of optimizing revision times by deriving and illustrating adaptive policies over the multi-period newsvendor problem.} We further develop a mixed-integer linear programming reformulation for the generic approach through scenario trees, and show that solving the resulting adaptive two-stage model is NP-Hard. Then, we focus our analyses on a specific problem structure that includes the capacity expansion planning problem under uncertainty. We derive the value of the proposed approach in comparison to {\color{black}the} two-stage and multi-stage stochastic programming models in terms of the revision times and problem parameters over the expansion planning problem. We also provide analytical analyses on these relative values with respect to the cost and demand parameters of the problem structure studied. 
\textcolor{black}{Furthermore, we propose a two-phase solution procedure by first determining the revision times and then solving the resulting program under these revision decisions.
For the capacity expansion planning problem, these algorithms select the revision times by minimizing the upper bound on the objective gap of the adaptive two-stage approach in comparison to the multi-stage and maximizing the lower bound of the corresponding gap with the two-stage models by benefiting from our analytical analyses. We present an additional solution algorithm based on the relaxation of the adaptive two-stage model to obtain the revision times, for which we develop a constant  approximation guarantee.
In order to illustrate our results, we study a generation capacity expansion planning problem over a multi-period planning horizon. 
Our extensive computational study highlights significant gains of the adaptive two-stage approach in objective in comparison to {\color{black}the} two-stage model, demonstrating the importance of optimizing revision decisions, with solution qualities that are almost similar to the multi-stage models.} 
We show that this relative gain even increases in scenario trees with higher variability, and  provide a computational study illustrating the performance and convergence of the proposed solution algorithms. Finally, we present the practical implications of utilizing the adaptive two-stage approach by examining sample generation capacity expansion plans. 

As a future research direction, our approach can be extended to various problem settings where limited flexibility is a requisite by determining the best time to revise decisions. 
{\color{black}Furthermore, tailored decomposition algorithms and decision rule approaches can be developed to address the solution of large-scale adaptive two-stage stochastic programs.} 
As an additional remark, the generic adaptive two-stage stochastic model 
may be investigated 
under stagewise independepence assumption of the random variables by exploring potential dynamic programming reformulations, solution algorithms, and closed form optimal policies for specific problem structures. 



\bibliographystyle{informs2014} 
\bibliography{References} 






\begin{APPENDICES}
{
\small

\section{Adaptive Two-stage Stochastic Programming for the Newsvendor Problem}
\label{adaptiveNewsvendorProof}

To obtain the result discussed in Theorem \eqref{ATSNewsvendorProp}, we first focus on the problem \eqref{eq:ATS_DPFormulation2} where the decision maker observes the inventory at hand at the beginning of period $t^*$ to determine the order schedule of periods ${t^*}, \cdots, T$. 
\begin{prop} \label{staticNewsvendorProp}
For the problem \eqref{eq:ATS_DPFormulation2}, we obtain an optimal policy in the following form:
\begin{equation}
\widetilde{F}_{t^*,t} (s + X_{t^*,t}) = \frac{-c_t + c_{t+1} + b_t}{h_t + b_t}, \quad t = t^*, \cdots, T.
\end{equation}
\end{prop} 
\proof{Proof:} 
We first observe that the objective function of the problem \eqref{eq:ATS_DPFormulation2} is convex as it can be represented as sum of convex functions. 
Then, we rewrite the objective function in the following form:
\begin{align*} 
H(x) & = \sum_{t=t^*}^T (c_t x_t + h_t \int_{0}^{s + \sum_{t'=t^*}^t x_{t' }} (s + \sum_{t'=t^*}^t x_{t' } - z) \tilde{f}_{t^*,t}(z)dz \\
& \hspace{60mm} - b_t \int_{s + \sum_{t'=t^*}^t x_{t' }}^{\infty} (s + \sum_{t'=t^*}^t x_{t' } - z)\tilde{f}_{t^*,t}(z) dz),
\end{align*}
where $\tilde{f}_{i,j}$ is the convolution probability density function corresponding to $\sum_{t=i}^j d_t$. We then take derivative of this expression with respect to the order quantities of the periods $t^*, \cdots, T$ as follows:
\begin{align*}
& \frac{\partial H(x)}{\partial x_{t^*}} = c_{t^*} + h_{t^*} \widetilde{F}_{t^*,t^*} (s+ x_t^*) - b_{t^*} (1 - \widetilde{F}_{t^*,t^*} (s + x_t^*)) + \cdots \\
& \hspace{70mm} + h_T \widetilde{F}_{t^*,T} (s + \sum_{t=t^*}^T x_t) - b_T (1 - \widetilde{F}_{t^*,T} (s + \sum_{t=t^*}^T x_t) ) \\
& \qquad \vdots \notag \\ 
& \frac{\partial H(x)}{\partial x_{T-1}} = c_{T-1} + h_{T-1} \widetilde{F}_{t^*,T-1} (s + \sum_{t=t^*}^{T-1} x_t) - b_{T-1} (1 - \widetilde{F}_{t^*,T-1} (s + \sum_{t=t^*}^{T-1} x_t)) \\
& \hspace{70mm} + h_T \widetilde{F}_{t^*,T} (s + \sum_{t=t^*}^T x_t) - b_T (1 - \widetilde{F}_{t^*,T} (s + \sum_{t=t^*}^T x_t)) \\
& \frac{\partial H(x)}{\partial x_T} = c_T +  h_T \widetilde{F}_{t^*,T} (s + \sum_{t=t^*}^T x_t) - b_T (1 - \widetilde{F}_{t^*,T} (s + \sum_{t=t^*}^T x_t))
\end{align*}
\textcolor{black}{To account for the nonnegativity of the order quantities, we introduce nonnegative variables $\lambda_t$ for $t = t^*,\cdots,T$ 
by considering the conditions $\partial H(x)/\partial x_{t} = \lambda_t$, $\lambda_t \> x_t = 0$ for $t = t^*,\cdots,T$. Then, we can conclude that any nonnegative $\{\lambda_t, x_t\}_{t = t^*}^T$ satisfying these conditions are optimal. Considering order up to policies on the order quantities $x_t$ in the form $x_{t^*} = \max\{X^*_{t^*,t^*},0\}$, and $x_t = \max\{X^*_{t^*,t} - \sum_{t' = t^*}^{t-1} x_{t'}, 0\}$ for $t = t^*+1, \cdots, T$, where $X^*_{t^*,t} = \widetilde{F}^{-1}_{t^*,t} (\frac{-c_t + c_{t+1} + b_t}{h_t + b_t}) - s$ for $t = t^*, \cdots, T$, 
we can obtain $\lambda_t$ for $t = t^*,\cdots,T$ satisfying the desired conditions and concluding the proof.} 
\Halmos\endproof

As a next step, we aim extending our result to the adaptive two-stage case. We observe that the objective function denoted in \eqref{eq:ATS_DPFormulation} is convex. Let the corresponding objective function be $G(x)$, and we take its derivative with respect to the order quantities of the periods $1, \cdots, t^* - 1$ as follows:
\begin{align*}
& \frac{\partial G(x)}{\partial x_{1}} = c_{1} + h_{1} \widetilde{F}_{1,1} (x_1) - b_{1} (1 - \widetilde{F}_{1,1} (x_1)) + \cdots \\
& \hspace{30mm} + h_{t^*-1} \widetilde{F}_{1, t^* - 1} (\sum_{t=1}^{t^*-1} x_t) - b_{t^*-1} (1 - \widetilde{F}_{1,t^* -1} (\sum_{t=1}^{t^*-1} x_t)) + \frac{\partial E[ Q_{t^*}(\sum_{t = 1}^{t^*-1} (x_t - d_t))]}{\partial x_1}\\
& \qquad \vdots \notag \\ 
& \frac{\partial G(x)}{\partial x_{t^*-1}} = c_{t^*-1} + h_{t^*-1} \widetilde{F}_{1,t^*-1} (\sum_{t=1}^{t^*-1} x_t) - b_{t^*-1} (1 - \widetilde{F}_{1,t^*-1} (\sum_{t=t^*}^{t^*-1} x_t)) + \frac{\partial E[Q_{t^*}(\sum_{t = 1}^{t^*-1} (x_t - d_t))]}{\partial x_{t^*-1}}
\end{align*}

To derive $\frac{\partial Q_{t^*}(s)}{\partial x_i}$ for $i = 1, \cdots, t^*-1$, we use the chain rule as $\frac{\partial Q_{t^*}(s)}{\partial x_i} = \frac{\partial Q_{t^*}(s)}{\partial s} . \frac{\partial s}{\partial x_i}$, where $s = \sum_{t = 1}^{t^*-1} (x_t - d_t)$. For that purpose, we need to identify the derivative of $Q_{t^*}(s)$ with respect to $s$.
\begin{lem}
$\displaystyle \frac{\partial Q_{t^*}(s)}{\partial s} = - c_{t^*}$.
\end{lem}
\proof{Proof:} 
We let the optimal solution of the problem \eqref{eq:ATS_DPFormulation2} as $(x^*_{t^*}, \cdots, x^*_T)$. Then, we can represent $\displaystyle \frac{\partial Q_{t^*}(s)}{\partial s}$ as follows:
\begin{align*}
& \frac{\partial Q_{t^*}(s)}{\partial s} = h_{t^*} \widetilde{F}_{t^*,t^*} (s+ x^*_{t^*}) - b_{t^*} (1 - \widetilde{F}_{t^*,t^*} (s + x^*_{t^*})) + \cdots \\
& \hspace{70mm} + h_T \widetilde{F}_{t^*,T} (s + \sum_{t=t^*}^T x^*_{t}) - b_T (1 - \widetilde{F}_{t^*,T} (s + \sum_{t=t^*}^T x^*_{t}) )
\end{align*}
Using the optimality conditions derived in the proof of Proposition \ref{staticNewsvendorProp}, we can prove that $\displaystyle \frac{\partial Q_{t^*}(s)}{\partial s} = - c_{t^*}$. 
\Halmos\endproof
\textcolor{black}{Combining the above with the nonnegativity conditions on $x_t$ by introducing the nonnegative variables $\lambda_t$ for $t = 1, \cdots, t^{*}-1$ with the relationships $\partial G(x)/\partial x_{t} = \lambda_t$, $\lambda_t x_t = 0$, for $t = 1, \cdots, t^{*}-1$, we can construct nonnegative $\{\lambda_t, x_t\}_{t = 1}^{t^{*}-1}$ pairs satisfying these conditions. To this end, we define the order up to policies in the form $x_1 = X^*_{1,1}$, and $x_t = \max\{X^*_{1,t} - \sum_{t' = 1}^{t-1} x_{t'}, 0\}$ 
for $t = 2, \cdots, t^*-1$ where $X^*_{1,t} = \widetilde{F}^{-1}_{1,t}(\frac{-c_t + c_{t+1} + b_t}{h_t + b_t})$ for $t = 1, \cdots, t^*-1$ for obtaining $\lambda_t$ variables satisfying the conditions. Thus, we obtain the desired result discussed in Theorem \eqref{ATSNewsvendorProp}.} 

Using Theorem \eqref{ATSNewsvendorProp}, we can show that the adaptive two-stage solution follows an order up to policy. 
Let $\{X^*_{1,t}\}_{t = 1}^T$ be the cumulative order up to quantities. {\color{black}To compute these values, we first compute $X^*_{1,t}$ for $t = 1, \cdots, t^*-1$, obtain $s_{t^*}$, and then compute $X^*_{1,t}$ for $t = t^*, \cdots, T$. Spefically, } 
$X^*_{1,t} = \widetilde{F}^{-1}_{1,t}(\frac{-c_t + c_{t+1} + b_t}{h_t + b_t})$ for $t = 1, \cdots, t^*-1$, and $X^*_{t^*,t} = \widetilde{F}^{-1}_{t^*,t} (\frac{-c_t + c_{t+1} + b_t}{h_t + b_t}) - s_{t^*}$ for $t = t^*, \cdots, T$. Next, we derive the order amount of each period as follows: $x^*_1 = X^*_{1,1}$, and $x^*_t = \max\{X^*_{1,t} - \sum_{t' = 1}^{t-1} x^*_{t'}, 0\}$ 
for $t = 2, \cdots, t^*-1$. At time $t^*$, we observe the cumulative net inventory of that period, $s_{t^*}$ as $\sum_{t' = 1}^{t^{*}-1} x^*_{t'} - D_{1,t^*-1}$. Then, we derive the remaining ordering policy as $x^*_{t^*} = \max\{X^*_{t^*,t^*},0\}$, and $x^*_t = \max\{X^*_{t^*,t} - \sum_{t' = t^*}^{t-1} x^*_{t'}, 0\}$ for $t = t^*+1, \cdots, T$. We note than when $t^*$ is set to 1, and $s_{t^*}$ represents the initial inventory, then the adaptive approach converts into a static setting where the decision maker determines the order schedules until the end of the planning horizon ahead of the planning. 

\section{Technical Proofs}
\label{appendixOmittedProofs}

\subsection{Proof of Theorem \ref{thmNPHard}}
\label{thmNPHardProof}

To prove the theorem, it suffices to show that the feasibility problem associated with the adaptive two-stage problem \eqref{eq:adaptiveTwostageFormulationNonlinear} is NP-Complete. For this purpose, we define the feasible region constructed in \eqref{eq:adaptiveTwostageFormulationNonlinearInstance}, and refer to this feasibility problem as P. Next, we demonstrate that the subset sum problem, which is known to be NP-Complete \citep{Garey1990}, can be reduced to  P 
in polynomial time. We specify the subset sum problem as follows: Given the non-negative integers $w'_1, w'_2, \cdots, w'_{N'}$, and $W'$, does there exist a subset $S \subseteq \{1, \cdots, N' \}$ such that $\sum_{i \in S} w'_i = W'$? 
\begin{subequations}  \label{eq:adaptiveTwostageFormulationNonlinearInstance}
\begin{alignat}{2}
\min_{\alpha, \beta, t^*} \quad & 0  \\
\text{s.t.} \quad & \alpha_{i2} - \beta_{i2} = \alpha_{i1} \quad & i = 1, \cdots, N \label{eq:constrList1}\\
& \alpha_{i3} - \beta_{i3} = -\alpha_{i1} \quad & i = 1, \cdots, N \label{eq:constrList1a}\\
& \alpha_{i4} - \beta_{i4} = \alpha_{i1} \quad & i = 1, \cdots, N \label{eq:constrList1b}\\
& \alpha_{i5} - \beta_{i5} = 2 - \alpha_{i1} \quad & i = 1, \cdots, N \\
& \alpha_{i6} - \beta_{i6} = \alpha_{i1} \quad & i = 1, \cdots, N \\
& \alpha_{i7} - \beta_{i7} = 2 - \alpha_{i1} \quad & i = 1, \cdots, N \label{eq:constrList2} \\
& \sum_{i = 1}^N w_i \alpha_{i1} = W & \label{eq:constrKnapsack} \\
& \alpha_{im} = \alpha_{in},  \beta_{im} = \beta_{in} \quad & \forall m, n \in S_t, \quad t < t_i^*, i = 1, \cdots, N \label{eq:constrNAC1} \\
& \alpha_{im} = \alpha_{in},  \beta_{im} = \beta_{in} \quad & \forall m, n \in S_t \cap {\mathcal{T}}(j), \quad j \in S_{t_i^*}, \quad t \geq t_i^*, i = 1, \cdots, N \label{eq:constrNAC2} \\
&  t_i^* \in \{1, 2, 3\}  & i = 1, \cdots, N. \label{eq:constrtStarRelationship}
\end{alignat}
\end{subequations}
Clearly, the problem \eqref{eq:adaptiveTwostageFormulationNonlinearInstance} is an instance of the adaptive two-stage problem. Specifically, we let $I = N$, $J = 0$, and $x_{in} =\begin{bmatrix} \alpha_{in} \\ \beta_{in} \end{bmatrix}$. We consider the scenario tree $\mathcal{T}$ with $T = 3$ stages 
as depicted in Figure \ref{NPHardScenarioTree}. We let $\mathcal{T} = \{1, \cdots, 7\}$ be the set of nodes in ascending order, where node 1 represents the root node. We also let $N = N'$, $W = W'$, and $w_i = w'_i$ for all $i = 1, \cdots, N$. 
We note that constraints \eqref{eq:constrList1} - \eqref{eq:constrList2} correspond to constraint \eqref{eq:compactConstraint1}, and constraint \eqref{eq:constrKnapsack} represent the constraint \eqref{eq:compactConstraint2}. Similarly, constraints \eqref{eq:constrNAC1}, \eqref{eq:constrNAC2} and \eqref{eq:constrtStarRelationship} refer to the constraints \eqref{eq:constrBeforeRevision1Nonlinear}, \eqref{eq:constrAfterRevision1Nonlinear} and \eqref{eq:tStarConstraintNonlinear}, respectively.

\begin{lem} \label{LemmaNPHard}
The following holds for the problem \eqref{eq:adaptiveTwostageFormulationNonlinearInstance}:
\begin{enumerate}
\item If $t_i^* = 1$ for any $i = 1, \cdots, N$, then the problem \eqref{eq:adaptiveTwostageFormulationNonlinearInstance} is infeasible. 
\item If $t_i^* = 2$, then $\alpha_{i1} = 1$ for all $i = 1, \cdots, N$.
\item If $t_i^* = 3$, then $\alpha_{i1} = 0$ for all $i = 1, \cdots, N$.
\end{enumerate}
\end{lem}

The proof of the Lemma \ref{LemmaNPHard} follows from the construction of the program \eqref{eq:adaptiveTwostageFormulationNonlinearInstance}. Specifically, when $t_i^* = 1$ for any $i = 1, \cdots, N$, then due to constraints \eqref{eq:constrNAC1} and  \eqref{eq:constrNAC2}, $\alpha_{i1} = - \alpha_{i1} = 2 - \alpha_{i1}$ resulting in infeasibility. When $t_i^* = 2$, then we have $\alpha_{i1} = 2 - \alpha_{i1}$ due to constraints \eqref{eq:constrList1b} and \eqref{eq:constrList2}, resulting in $\alpha_{i1} = 1$. Similarly, when $t_i^* = 3$, we have $\alpha_{i1} = 2 - \alpha_{i1}$ due to constraints \eqref{eq:constrList1} and \eqref{eq:constrList1a} making $\alpha_{i0} = 0$.

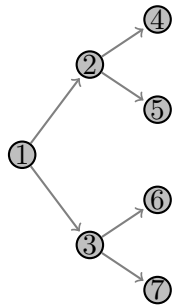
\begin{figure}[h]
\caption{Scenario tree for the problem \eqref{eq:adaptiveTwostageFormulationNonlinearInstance}.}
\label{NPHardScenarioTree}
\centering
\begin{tikzpicture}[shorten >=1pt,->,draw=black!50, node distance=\layersep, thick, scale=0.6]
    \tikzstyle{every pin edge}=[<-,shorten <=1pt]
    \tikzstyle{neuron}=[circle,fill=black!75,minimum size=10pt,draw=black,inner sep=0pt]
    \tikzstyle{input neuron}=[neuron, fill=black!25];
    \tikzstyle{output neuron}=[neuron, fill=white!00, draw=white];
    \tikzstyle{hidden neuron}=[neuron, fill=blue!50];
    \tikzstyle{annot} = [text width=4em, text centered]
	
	\node[input neuron] (I0) at (0,0) {1};
	\node[input neuron] (I1) at (1.5,2) {2};
	\node[input neuron] (I2) at (1.5,-2) {3};
	\node[input neuron] (I3) at (3,3) {4};
	\node[input neuron] (I4) at (3,1) {5};
	\node[input neuron] (I5) at (3,-1) {6};
	\node[input neuron] (I6) at (3,-3) {7};
	
            \path (I0) edge (I1);
            \path (I0) edge (I2);
            \path (I1) edge (I3);
            \path (I1) edge (I4);
            \path (I2) edge (I5);
            \path (I2) edge (I6);
    
\end{tikzpicture}
\end{figure}

We first show that given a solution $S$ to the subset-sum problem, we can construct a feasible solution for problem P. 
Specifically, we set $\alpha_{i1} = 1$ if $i \in S$, and 0 otherwise for all $i = 1, \cdots, N$. If $\alpha_{i1} = 1$, then $t_i = 2$ resulting in $\alpha_{i2} = 1$, $\beta_{i2} = 0$, $\alpha_{i3} = 0$, $\beta_{i3} = 1$, $\beta_{i4} = \beta_{i5} = \beta_{i6} = \beta_{i7} = 0$, $\alpha_{i4} = \alpha_{i5} = \alpha_{i6} = \alpha_{i7} = 1$ for all $i = 1, \cdots, N$. If  $\alpha_{i1} = 0$, then $t_i = 3$ resulting in $\alpha_{i2} = \beta_{i2} = \alpha_{i3} = \beta_{i3} = 0$, $\alpha_{i4} = \beta_{i4} = 0$, $\alpha_{i5} = 2$, $\beta_{i5} = 0$, $\alpha_{i6} = 0$, $\beta_{i6} = 0$, $\alpha_{i7} = 2$, $\beta_{i7} = 0$. Secondly, given a feasible solution $(\alpha, \beta, t^*)$ to the problem P, we can construct a feasible solution $S$ for the subset-sum. We note that in a feasible solution of P, we have either $t^*_i$ = 2 or $t_i^* = 3$ for $i = 1, \cdots, N$, as noted in Lemma \ref{LemmaNPHard}. To construct a feasible solution, we select $S = \{ i: t^*_i = 2\}$. As $(\alpha, \beta, t^*)$ is a feasible solution of P, we satisfy the condition $\sum_{i \in S} w'_i = W'$. 
Combining the above, we prove that 
$P$ is NP-Complete, completing the proof.

\subsection{Proof of Proposition \ref{TUAdaptiveProp}}
\label{TUAdaptivePropProof}

We observe that the coefficient matrix of \eqref{eq:CapacityExpansionAdaptiveReformulationCondensed} is composed of 0 and 1 values. Additionally, number of 1's in each column~$j$ is equal to $|{\hat{\mathcal{T}}}_i (t_i^*)(j)|$, where ${\hat{\mathcal{T}}}_i (t_i^*)(j)$ represents the subtree rooted at node~$j$. More specifically, an entry corresponding to $i^{th}$ row and $j^{th}$ column of the coefficient matrix is 1 if and only if node $i \in {\hat{\mathcal{T}}}_i (t_i^*)(j)$. Consequently, we can rearrange the coefficient matrix to obtain an interval matrix. This demonstrates that the desired matrix is totally unimodular. 

\subsection{Proof of Proposition \ref{PropCapExpBoundAdaptive}}
\label{PropCapExpBoundAdaptiveProof}

The decision variables $\{x_{n}\}_{\{n \in {\mathcal{T}}\}}$ can be categorized into two groups according to the revision time $t^*$. 
First of all, for $t < t^*$, the solutions $\{x_{m}: m \in S_t\}$ have the same value. Thus, we denote variables until period $t^*$ as $\hat{x}_{t}$ for all $t = 1, \cdots, t^*- 1$. Secondly, for $t \geq t^*$, we have the same solutions for $\{x_{m}: m \in S_{t} \cap {\mathcal{T}}(n)\}$ for all $n \in S_{t^*}$. We refer to these variables as $\hat{x}_{nt}$ for any $n \in S_{t^*}$ and $t \geq t^*$.

In order to find an upper bound on $v^R (t^*)$, we construct a feasible solution where $\hat{x}_{1} = \delta_{1}$, $\hat{x}_{t} = \max_{m \in {\mathcal{T}}(1, t)} \{\delta_{m} \} - \max_{m \in {\mathcal{T}}(1, t - 1)} \{\delta_{m} \}$ for $2 \leq t \leq t^* - 1$, and $\hat{x}_{nt} = \max \{\delta_{m}: m \in {\mathcal{T}}(1, t^*-1) \cup {\mathcal{T}}(n, t)\} - \max \{\delta_{m}: m \in {\mathcal{T}}(1, t^*-1) \cup {\mathcal{T}}(n, t-1)\}$ for $t \geq t^*$. Using the relationship $v^R (t^*) \leq \sum_{n \in {\mathcal{T}}} p_n  a_{n} x_{n}$ for any feasible $\{x_{n}\}_{\{n \in {\mathcal{T}}\}}$, we can obtain the following: 
\begin{align*}
v^R (t^*) & \leq \sum_{t=1}^{t^*-1} \sum_{n \in S_t} p_n  a_{n} \hat{x}_{t} + \sum_{t=t^*}^{T} \sum_{n \in S_{t^*}} p_n  a_{n} \hat{x}_{nt} \\
& \leq \bar{a}^- (t^*) \sum_{t=1}^{t^*-1} \hat{x}_{t} + \bar{a}^+ (t^*) \sum_{n \in S_{t^*}} p_n \sum_{t=t^*}^{T} \hat{x}_{nt} \\
& = \bar{a}^- (t^*) \sum_{t=1}^{t^*-1} ( \max_{m \in {\mathcal{T}}(1, t)} \{\delta_{m} \} - \max_{m \in {\mathcal{T}}(1, t - 1)} \{\delta_{m} \} ) \\
& \qquad +   \bar{a}^+ (t^*) \sum_{n \in S_{t^*}} p_n \sum_{t=t^*}^{T} (\max_{m \in {\mathcal{T}}(1, t^*-1) \cup {\mathcal{T}}(n, t)}\{\delta_{m}\} - \max_{m \in {\mathcal{T}}(1, t^*-1) \cup {\mathcal{T}}(n, t-1)} \{\delta_{m}\}) \\
& = \bar{a}^- (t^*) \max_{m \in {\mathcal{T}}(1, t^*-1)} \{\delta_{m}\} + \> \bar{a}^+ (t^*) \sum_{n \in S_{t^*}} p_n (\max_{m \in {\mathcal{T}}(1, t^*-1) \cup {\mathcal{T}}(n)} \{\delta_{m}\} - \max_{m \in {\mathcal{T}}(1, t^*-1)} \{\delta_{m}\}) \\
& = \bar{a}^- (t^*) \max_{m \in {\mathcal{T}}(1, t^*-1)} \{\delta_{m}\} + \> \bar{a}^+ (t^*) \sum_{n \in S_{t^*}} p_n \max_{m \in {\mathcal{T}}(1, t^*-1) \cup {\mathcal{T}}(n)} \{\delta_{m}\} - \> \bar{a}^+ (t^*) \max_{m \in {\mathcal{T}}(1, t^*-1)} \{\delta_{m}\}  \\
& = (\bar{a}^- (t^*) - \bar{a}^+(t^*)) \delta^-(t^*) + \bar{a}^+(t^*) \delta^+(t^*).
\end{align*}

Let $\hat{x}_{t}^*$ for $t < t^*$ and $\hat{x}_{nt}^*$, for $t \geq t^*$ and $n \in S_{t^*}$ denote the optimal solution for the problem \eqref{eq:CapacityExpansionAdaptiveReformulation}. Then, we can derive the following:
\begin{align*}
v^R (t^*) & = \sum_{t=1}^{t^*-1} \sum_{n \in S_t} p_n  a_{n} \hat{x}_{t}^* + \sum_{t=t^*}^{T} \sum_{n \in S_{t^*}} p_n  a_{n} \hat{x}_{nt}^* \\
& \geq \underbar{$a$}^- (t^*) \sum_{t=1}^{t^*-1} \hat{x}_{t}^* + \underbar{$a$}^+ (t^*) \sum_{n \in S_{t^*}} p_n \sum_{t=t^*}^{T} \hat{x}_{nt}^* \\
& \geq \underbar{$a$}^- (t^*) \sum_{t=1}^{t^*-1} \hat{x}_{t}^* + \underbar{$a$}^+ (t^*) \sum_{n \in S_{t^*}} p_n (\max_{m \in {\mathcal{T}}(1, t^* -  1) \cup {\mathcal{T}}(n)} \{\delta_{m} \} - \sum_{t = 1}^{t^*-1} \hat{x}_{t}^*) \\ 
& = (\underbar{$a$}^- (t^*) - \underbar{$a$}^+ (t^*) ) \sum_{t = 1}^{t^*-1} \hat{x}_{t}^* + \underbar{$a$}^+ (t^*) \sum_{n \in S_{t^*}} p_n (\max_{m \in {\mathcal{T}}(1, t^* -  1) \cup {\mathcal{T}}(n)} \{\delta_{m} \}) \\
& \geq (\underbar{$a$}^- (t^*) - \underbar{$a$}^+ (t^*) ) \delta^- (t^*) + \underbar{$a$}^+ (t^*)   \delta^+ (t^*).
\end{align*}

The second and fourth inequalities follow from constraint \eqref{eq:SubstructureConstr1}. 
In particular, any feasible solution $\hat{x}_{t}$ and $\hat{x}_{nt}$ for any $n \in S_{t^*}$ should satisfy $\sum_{t = 1}^{t^*-1} \hat{x}_{t} + \sum_{t=t^*}^T \hat{x}_{nt} \geq \max_{m \in {\mathcal{T}}(1, t^* -  1) \cup {\mathcal{T}}(n)} \{\delta_{m} \}$. Additionally, $\sum_{t = 1}^{t^*-1} \hat{x}_{t} \geq \max_{m \in {\mathcal{T}}(1, t^* -  1)} \{\delta_{m} \} = \delta^- (t^*)$. 

Combining above, we derive upper and lower bounds to $v^R (t^*)$.

\subsection{Proof of Proposition \ref{demandSensitivityPropGeneric}}
\label{demandSensitivityPropGenericProof}

If $\hat t^D = 1$, then $\delta^+(t) = \delta^*$ for any $t \in \{1, \cdots, T\}$. Otherwise, $\delta^+(1) = \delta^*$ as ${\mathcal{T}}(1)$ corresponds to the full scenario tree including the maximum demand value. For any node $n$ in stages $t > \hat t^D$, we have $\max_{m \in {\mathcal{T}}(1, t -  1) \cup {\mathcal{T}}(n)} \{\delta_{m} \} = \delta^*$. Hence, we obtain $\delta^+(t) = \delta^*$ for $t > \hat t^D$. For an intermediate stage $t$ between 1 and $\hat t^D$, we have the relationship $\delta^+(t) \leq \delta^*$. Combining the above, under a general demand structure for $\delta_{n}$ when $\hat t^D > 1$, we conclude that the minimizer of $\delta^+(t)$ is in $\{2, \cdots, \hat t^D\}$.

\subsection{Proof of Proposition \ref{demandSensitivityPropSpecific}}
\label{demandSensitivityPropSpecificProof}
We observe that 
for the cases $t < \hat t^D$ and $t > \hat t^D$, we have $\delta^+(t) = \delta^*$. 
For $t = t^D$, we obtain the relationship $\max_{m \in {\mathcal{T}}(1, t - 1) \cup {\mathcal{T}}(n)} \{\delta_{m} \} < \delta^*$ for some $n \in S_{\hat t^D}$. Hence, we have $\delta^+(\hat t^D) < \delta^*$ making $\hat t^D$ the minimizer of $\delta^+(t)$.

\subsection{Proof of Proposition \ref{capExpansionBoundProp1}}
\label{capExpansionBoundProp1Proof}

We first observe that $V^{TS} \geq \sum_{n \in {\mathcal{T}}} p_n b_n^\top y_n^{TLP} + \sum_{i=1}^{I} v_i^T$. 
Next, we note that $\{y_n^{TLP}\}_{n \in {\mathcal{T}}}$ is a feasible solution to the adaptive two-stage problem under $t^*$. Thus, $V^{ATS}(t^*) \leq \sum_{n \in {\mathcal{T}}} p_n b_n^\top y_n^{TLP} + \sum_{i=1}^{I} o_i^R (t_i^*)$, where  
\begin{equation*}
o_i^R (t_i^*) = \min \{\sum_{n \in {\hat{\mathcal{T}}}_i (t_i^*)} \hat p_n  \hat a_{in} x_{in}: \quad \sum_{m \in \hat{P}(n)} x_{im} \geq \hat \delta_{in}, x_{in} \in \mathbb{Z}_{+} \>  \forall n \in {\hat{\mathcal{T}}}_i (t_i^*) \}.
\end{equation*} 

Additionally, we can represent $v_i^R (t_i^*)$ as follows
\begin{align*}
v_i^R (t_i^*) & = \min \{\sum_{n \in {\hat{\mathcal{T}}}_i (t_i^*)} \hat p_n  \hat a_{in} x_{in}: \quad \sum_{m \in \hat{P}(n)} x_{im} \geq \hat \delta_{in}, x_{in} \in \mathbb{R}_{+} \>  \forall n \in {\hat{\mathcal{T}}}_i (t_i^*) \} \\
& = \max \{\sum_{n \in {\hat{\mathcal{T}}}_i (t_i^*)} \hat \delta_{in} \pi_{in}: \quad \sum_{m \in {\hat{\mathcal{T}}}_i (t_i^*)} \pi_{im} \leq \hat p_n \hat a_{in}, \pi_{in} \in \mathbb{R}_{+} \>  \forall n \in {\hat{\mathcal{T}}}_i (t_i^*) \}.
\end{align*}

Using Proposition \ref{TUAdaptiveProp} and linear programming duality, we can reexpress $o_i^R (t_i^*)$ as 
\begin{align*}
o_i^R (t_i^*) & = \min \{\sum_{n \in {\hat{\mathcal{T}}}_i (t_i^*)} \hat p_n  \hat a_{in} x_{in}: \quad \sum_{m \in \hat{P}(n)} x_{im} \geq \lceil \hat \delta_{in} \rceil, x_{in} \in \mathbb{R}_{+} \>  \forall n \in {\hat{\mathcal{T}}}_i (t_i^*) \} \\
& = \max \{\sum_{n \in {\hat{\mathcal{T}}}_i (t_i^*)} (\hat \delta_{in} + (\lceil \hat \delta_{in} \rceil - \hat \delta_{in})) \pi_{in}: \quad \sum_{m \in {\hat{\mathcal{T}}}_i (t_i^*)} \pi_{im} \leq \hat p_n \hat a_{in}, \pi_{in} \in \mathbb{R}_{+} \>  \forall n \in {\hat{\mathcal{T}}}_i (t_i^*) \} \\
& \leq v_i^R (t_i^*) + \max_{n \in {\hat{\mathcal{T}}}_i (t_i^*)} \{\lceil \hat \delta_{in} \rceil - \hat \delta_{in}\} \> \max \{\sum_{n \in {\hat{\mathcal{T}}}_i (t_i^*)} \pi_{in}: \quad \sum_{m \in {\hat{\mathcal{T}}}_i (t_i^*)} \pi_{im} \leq \hat p_n \hat a_{in}, \pi_{in} \in \mathbb{R}_{+} \>  \forall n \in {\hat{\mathcal{T}}}_i (t_i^*) \} \\
& = v_i^R (t_i^*) + \max_{n \in {\hat{\mathcal{T}}}_i (t_i^*)} \{\lceil \hat \delta_{in} \rceil - \hat \delta_{in}\} \> \min \{\sum_{n \in {\hat{\mathcal{T}}}_i (t_i^*)} \hat \hat p_n  \hat a_{in} x_{in}: \quad\sum_{m \in \hat{P}(n)} x_{im} \geq 1, x_{in} \in \mathbb{R}_{+} \>  \forall n \in {\hat{\mathcal{T}}}_i (t_i^*) \} \\
& = v_i^R (t_i^*) + \max_{n \in {\hat{\mathcal{T}}}_i (t_i^*)} \{\lceil \hat \delta_{in} \rceil - \hat \delta_{in}\} \> a_{i1}.
\end{align*}
Here, the last equality follows from the fact that $x_{i1} = 1$ and $x_{in} = 0$ for $n \in {\hat{\mathcal{T}}}_i (t_i^*) \setminus \{1\}$ is an optimal solution of the resulting single-resource problem. Combining the above, we demonstrate the desired result. 

\subsection{Proof of Proposition \ref{capExpansionBoundProp2}}
\label{capExpansionBoundProp2Proof}

We first observe that $V^{MS} \geq \sum_{n \in {\mathcal{T}}} p_n b_n^\top y_n^{MLP} + \sum_{i=1}^{I} v_i^M$. Next, we have $V^{ATS}(t^*) \leq \sum_{n \in {\mathcal{T}}} p_n b_n^\top y_n^{MLP} + \sum_{i=1}^{I} o_i^R (t_i^*)$ as before. 
Using the same techniques as in the proof of Proposition~\ref{capExpansionBoundProp1}, we obtain the desired result. 

\subsection{Proof of Proposition \ref{approxGuarenteeProposition}}
\label{capExpansionApproximationProof}

We denote the objective value of the adaptive two-stage version of the stochastic capacity expansion planning problem \eqref{eq:CapacityExpansionMultiStage} under a given $(x,y,t)$ decision as \[
f(x,y,t) := \sum_{i=1}^I \sum_{n \in {\mathcal{T}}_i (t_i)} \hat p_n \hat a_{in} x_{in} + \sum_{n \in {\mathcal{T}}} p_n b_n^\top y_n,
\]
where ${\mathcal{T}}_i (t_i)$ corresponds to the compressed tree under the revision decision $t_i$, 
as discussed in formulation~\eqref{eq:CapacityExpansionAdaptiveReformulationCondensed}. 

We represent the solution corresponding to ATS-Relax algorithm as $(x^{ATS-Relax},$ $y^{ATS-Relax},$  $t^{ATS-Relax})$, and the solution of the true adaptive two-stage program as $(x^{ATS}, y^{ATS}, t^{ATS})$. Consequently, we define the bound between the two approaches as
\begin{align*}
V^{ATS-Relax} - V^{ATS} = & f (x^{ATS-Relax}, y^{ATS-Relax}, t^{ATS-Relax}) - f(x^{ATS}, y^{ATS}, t^{ATS}) \\
\leq & f(x^{ATS-Relax}, y^{ATS-Relax}, t^{ATS-Relax}) - f(x^{LP}, y^{LP}, t^{LP}) \\
= & f(x^{ATS-Relax}, y^{ATS-Relax}, t^{ATS-Relax}) - f(x^{ATS-Relax}, y^{LP}, t^{LP}) \\
& + f(x^{ATS-Relax}, y^{LP}, t^{LP}) - f(x^{LP}, y^{LP}, t^{LP}) \\
\leq & f(x^{ATS-Relax}, y^{LP}, t^{LP}) - f(x^{LP}, y^{LP}, t^{LP}),
\end{align*}
where $(x^{LP}, y^{LP}, t^{LP})$ represents the solution of the relaxation of the true adaptive two-stage program when $x$ decisions are relaxed to be continuous. We note that $t^{LP} = t^{ATS-Relax}$ by the definition of Algorithm \ref{alg:AdaptiveApproxAlg3}. Next, we can state the following 
\begin{equation*}
f(x^{ATS-Relax}, y^{LP}, t^{LP}) - f(x^{LP}, y^{LP}, t^{LP}) = \sum_{i=1}^I \sum_{n \in {\mathcal{T}}_i(t_i^{LP})} \hat p_n \hat a_{in} (x_{in}^{ATS-Relax} - x_{in}^{LP}),
\end{equation*}
Here, $\hat p_n$ and $\hat a_n$ are computed specifically for the compressed tree under the revision decision $t^{LP}$. Using the analysis in the proof of Theorem 6 in \citet{Huang2009}, the above expression reduces to
\begin{equation*}
f(x^{ATS-Relax}, y^{LP}, t^{LP}) - f(x^{LP}, y^{LP}, t^{LP}) = \sum_{i = 1}^I (\max_{n \in {\mathcal{T}}_i (t_i^{ATS-Relax})} \{(\lceil {\color{black}\hat \delta_{in} \rceil - \hat \delta_{in}})\} a_{i1}).
\end{equation*}
Combining the above, we obtain the desired result. 

\section{Illustrative Instance}
\label{appendixInstanceData}

In this section, we provide the details of the instance studied in Section \ref{senstivityVatsSingleResource}. The cost parameter $a_n = 1$ for all $n \in {\mathcal{T}}$. The demand parameter $\{\delta_n\}_{n \in {\mathcal{T}}}$ is randomly generated from the distribution $N(\mu,\sigma^2)$, where $\mu$ = 30 and $\sigma$ = 5. We consider a scenario tree with 5 stages, and the resulting values for $\{\delta_n\}_{n \in {\mathcal{T}}}$ are presented in Figure \ref{IllustrativeInstanceDeltaTree}.

\begin{figure}[H]
\caption{Demand values for the illustrative instance.}
\label{IllustrativeInstanceDeltaTree}
    \centering
\begin{tikzpicture}[shorten >=1pt,->,draw=black!50, node distance=\layersep, thick, scale=0.25]
    \tikzstyle{every pin edge}=[<-,shorten <=1pt]
    \tikzstyle{neuron}=[minimum size=5pt,inner sep=2pt, circle] 
    \tikzstyle{input neuron}=[neuron];
    \tikzstyle{output neuron}=[neuron, fill=white!00, draw=white];
    \tikzstyle{hidden neuron}=[neuron, fill=blue!50];
    \tikzstyle{annot} = [text width=4em, text centered]

	\node[input neuron] (I0) at (0,0) {27};
	\node[input neuron] (I1) at (7,4) {29};
	\node[input neuron] (I2) at (7,-4) {19};
	\node[input neuron] (I3) at (14,7.5) {38};
	\node[input neuron] (I4) at (14,2.5) {21};
	\node[input neuron] (I5) at (14,-2.5) {25};
	\node[input neuron] (I6) at (14,-7.5) {32};
	\node[input neuron] (I7) at (21,10.5) {23};
	\node[input neuron] (I8) at (21,7.5) {24};
	\node[input neuron] (I9) at (21,4.5) {25};
	\node[input neuron] (I10) at (21,1.5) {32};
	\node[input neuron] (I11) at (21,-1.5) {41};
	\node[input neuron] (I12) at (21,-4.5) {30};
	\node[input neuron] (I13) at (21,-7.5) {24};
	\node[input neuron] (I14) at (21,-10.5) {32};
	\node[input neuron] (I15) at (28,11.25) {27};
	\node[input neuron] (I16) at (28,9.75) {29};
	\node[input neuron] (I17) at (28,8.25) {35};
	\node[input neuron] (I18) at (28,6.75) {26};
	\node[input neuron] (I19) at (28,5.25) {30};
	\node[input neuron] (I20) at (28,3.75) {25};
	\node[input neuron] (I21) at (28,2.25) {29};
	\node[input neuron] (I22) at (28,0.75) {31};
	\node[input neuron] (I23) at (28,-0.75) {25};
	\node[input neuron] (I24) at (28,-2.25) {28};
	\node[input neuron] (I25) at (28,-3.75) {28};
	\node[input neuron] (I26) at (28,-5.25) {26};
	\node[input neuron] (I27) at (28,-6.75) {24};
	\node[input neuron] (I28) at (28,-8.25) {22};
	\node[input neuron] (I29) at (28,-9.75) {28};
	\node[input neuron] (I30) at (28,-11.25) {29};
	
            \path (I0) edge (I1);
            \path (I0) edge (I2);
            \path (I1) edge (I3);
            \path (I1) edge (I4);
            \path (I2) edge (I5);
            \path (I2) edge (I6);
            \path (I3) edge (I7);
            \path (I3) edge (I8);
            \path (I4) edge (I9);
            \path (I4) edge (I10);
            \path (I5) edge (I11);
            \path (I5) edge (I12);
            \path (I6) edge (I13);
            \path (I6) edge (I14);
            \path (I7) edge (I15);
            \path (I7) edge (I16);
            \path (I8) edge (I17);
            \path (I8) edge (I18);
            \path (I9) edge (I19);
            \path (I9) edge (I20);
            \path (I10) edge (I21);
            \path (I10) edge (I22);
            \path (I11) edge (I23);
            \path (I11) edge (I24);
            \path (I12) edge (I25);
            \path (I12) edge (I26);
            \path (I13) edge (I27);
            \path (I13) edge (I28);
            \path (I14) edge (I29);
            \path (I14) edge (I30);
    
\end{tikzpicture}
\end{figure}
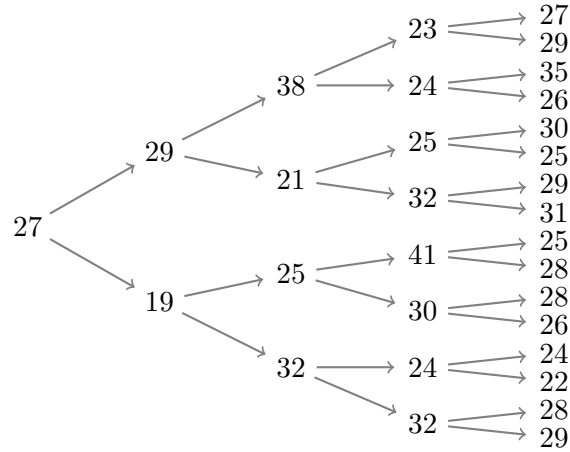

\section{Cost Sensitivity for the Single-Resource Problem}
\label{CostSensitivitySectionAppendix}

In this section, we consider the case where the demand parameters $\{\delta_{n}\}_{n \in {\mathcal{T}}}$ are equal to each other throughout the scenario tree, i.e., {\color{black}$\{\delta_{n}\}_{n \in {\mathcal{T}}} = \delta$} for the single-resource problem studied in Section \ref{senstivityVatsSingleResource}. 
In this case, the bounds \eqref{eq:subBoundAdjustable} and \eqref{eq:subBoundMultistage} reduce to the following 
{\color{black}
\begin{align}
(a_{*} - \bar{a}^- (t^*))\delta & \leq v^T - v^R (t^*) \leq (a^* -  \underbar{$a$}^- (t^*)) \delta, \label{eq:CostSensitiveBoundTS} \\
(\underbar{$a$}^- (t^*) - a^*)\delta & \leq v^R (t^*) - v^M \leq (\bar{a}^- (t^*) - a_{*}) \delta. \label{eq:CostSensitiveBoundMS}
\end{align}
}
We observe that the value of the adaptive formulation depends on the choice of the revision point as the bounds in \eqref{eq:CostSensitiveBoundTS} and \eqref{eq:CostSensitiveBoundMS} are functions of $t^*$.
\textcolor{black}{We note that the lower bounds are nonpositive for this setting, considering the definitions of the cost parameters.}
In order to select the best revision point for the adaptive two-stage approach, we aim \textcolor{black}{at} maximizing the lower bound on its objective difference between two-stage in \eqref{eq:CostSensitiveBoundTS} and \textcolor{black}{minimizing} the upper bound on the corresponding difference between multi-stage in \eqref{eq:CostSensitiveBoundMS}. Thus, this results in finding the revision point $2 \leq t^* \leq T$ that minimizes $\bar{a}^- (t^*) - a_{*}$. Since $a_{*}$ is not dependent to the choice of the revision point, it is omitted in the remainder of our analysis. We denote the best revision point in terms of the cost bounds as $t^{CB} := \argmin_{2 \leq t \leq T} \bar{a}^- (t)$. 

We consider several specific cases to illustrate the revision point decision $t^{CB}$ based on the cost values. If cost values are increasing in each stage of the scenario tree, then $\bar{a}^- (t)$ is monotonically increasing in $t$. Thus, the revision point needs to be as early as possible. On the other hand, if cost values are decreasing in each stage, then 
the value of the adaptive approach has the same value in all time periods as  $\bar{a}^- (t)$ is same for every $2 \leq t \leq T$.
Thus, revision time does not affect the analytical bounds for this setting. We can generalize this result as follows by letting the maximum cost value over the scenario tree as $\hat t^C = \min \{t \in \{1, \cdots, T\}: a_{n} = a^*, n \in S_t \}$.

\begin{prop} \label{costSensitiveProposition}
Under a general cost structure for $a_{n}$ with demand values {\color{black}$\{\delta_{n}\}_{n \in {\mathcal{T}}} = \delta$}, $\bar{a}^- (t)$ is monotonically nondecreasing until the period $\hat t^C$, and it remains constant afterwards. 
\end{prop}

The proof of Proposition \ref{costSensitiveProposition} is immediate by following the definition of $\bar{a}^- (t)$. This proposition shows that if the maximum cost value is in the root node of the scenario tree, then the value of the adaptive approach is not affected by the revision decision. Otherwise, the revision point needs to be selected as early as possible before observing the highest cost value of the scenario tree. 

\section{Scenario Tree Generation Algorithm}
\label{scenarioTreeSectionAppendix}

\begin{algorithm}[H]
\small
\caption{Scenario Tree Generation with $M$ branches}
\label{alg:scenarioTreeAlgorithm}
\begin{algorithmic}[1]
\STATE Obtain the demand values at the root node as $\{d_{k0}\}_{k \in \mathcal{K}}$.
\FORALL{$t = 2, \cdots, T$}
\STATE Let $\underline{\alpha}_t \leq \overline{\alpha}_t $ be the demand increase multiplier bounds.
\STATE Find the equisized multiplier interval values as $\{\alpha_t^j\}_{j = 0}^{M}$ where $\alpha_t^j = \underline{\alpha}_t + j (\overline{\alpha}_t -  \underline{\alpha}_t)/M$.
\FORALL{$k \in \mathcal{K}$}
\FORALL{node $n \in S_t$}
\STATE Find the order of the node in the stage as $j = n \> mod \> M$.
\STATE Generate the demand increase multiplier $\beta_t^j \sim U(\alpha_t^j, \alpha_t^{j+1})$.
\STATE Define $d_{k n} := \beta_t^j d_{k a(n)}$. 
\ENDFOR
\ENDFOR
\ENDFOR
\end{algorithmic}
\end{algorithm}

}
\end{APPENDICES}

\end{document}